\let\ORIlabel\label
\let\ORIrefstepcounter\refstepcounter
\AddToHook{package/hyperref/before}{%
  \let\label\ORIlabel
  \let\refstepcounter\ORIrefstepcounter
}
\documentclass[onefignum,onetabnum,final,dvipsnames]{siamart190516}
\usepackage[T1]{fontenc}
\usepackage[utf8]{inputenc}
\usepackage{graphicx}
\usepackage{dsfont}
\usepackage{amssymb}
\usepackage{amsfonts}
\usepackage{mathrsfs}
\usepackage{enumitem}
\usepackage{etoolbox}
\usepackage{subcaption}
\usepackage{booktabs}
\usepackage{tabularx}
\usepackage{mathtools}
\usepackage{tgpagella}
\usepackage{aligned-overset}
\usepackage{bm}
\usepackage{todonotes}
\usepackage{wasysym} 
\usepackage[most]{tcolorbox}
\usepackage{tikz}
\usepackage{bibentry}
\usepackage[dvipsnames]{xcolor}
\usepackage[normalem]{ulem}

\newcommand{\Hol}{\mathrm{Hol}}

\newcommand{\rev}[1]{#1}
\newcommand{\revsout}[1]{}
\usetikzlibrary{matrix}
\let\vec\oldvec

\let\epsilon\varepsilon

\let\phi\varphi

\DeclareMathOperator{\spn}{span}

\DeclareMathOperator{\vec}{\mathbf{vec}}
\DeclareMathOperator{\tvec}{\widetilde{\mathbf{vec}}}
\DeclareMathOperator{\matr}{\mathbf{matr}}

\DeclareMathOperator{\depth}{L}
\DeclareMathOperator{\size}{M}
\DeclareMathOperator{\Realiz}{R}

\newsiamremark{remark}{Remark}

\newcommand{\bvec}[1]{\bm{#1}}
\newcommand{\bmat}[1]{\mathbf{#1}}

\newcommand{\lnorm}{\left\vert\kern-0.25ex\left\vert\kern-0.25ex\left\vert}
\newcommand{\rnorm}{\right\vert\kern-0.25ex\right\vert\kern-0.25ex\right\vert}

\newcommand{\tC}{{\widetilde{C}}}

\newcommand{\bk}{{\bm{k}}}

\newcommand{\tu}{\widetilde{u}}

\newcommand{\cE}{\mathcal{E}}

\newcommand{\bu}{\bm{u}}
\newcommand{\bv}{\bm{v}}
\newcommand{\bx}{\bvec{x}}
\newcommand{\bb}{\bvec{b}}
\newcommand{\by}{\bvec{y}}

\newcommand{\bw}{\bm{w}}
\newcommand{\bbf}{\bm{f}}
\newcommand{\bzero}{{\bvec{0}}}

\newcommand{\is}{{i_1,\dots, i_d}}

\newcommand{\tphi}{\widetilde{\phi}}

\newcommand{\brange}[1]{\{1, \dots, #1\}}

\newcommand{\cqone}{c_{\mathrm{quad},1}}
\newcommand{\cqtwo}{c_{\mathrm{quad},2}}
\newcommand{\wkq}{w^{(q)}_k}
\newcommand{\wkppone}{w^{(p+1)}_k}
\newcommand{\xq}{\bvec{x}^{(q)}}
\newcommand{\xkq}{\bvec{x}^{(q)}_k}
\newcommand{\xkppone}{\bvec{x}^{(p+1)}_k}
\newcommand{\Cp}{C_{\mathrm{p}}}
\newcommand{\bp}{b_{\mathrm{p}}}
\newcommand{\alphap}{\alpha_{\mathrm{p}}}

\newcommand{\cA}{\mathcal{A}}
\newcommand{\genpar}{\theta}
\newcommand{\genparmin}{\genpar_{\mathrm{min}}}
\newcommand{\genset}{\Theta}
\newcommand{\cR}{\mathcal{R}}
\newcommand{\cP}{\mathcal{P}}
\newcommand{\cD}{\mathcal{D}}

\newcommand{\cU}{\mathcal{U}}
\newcommand{\cO}{\mathcal{O}}

\newcommand{\hbD}{\bmat{\widehat{D}}}
\newcommand{\hbM}{\bmat{\widehat{M}}}
\newcommand{\tbD}{\bmat{\widetilde{D}}}
\newcommand{\N}{\mathbb{N}}
\newcommand{\Z}{\mathbb{Z}}
\newcommand{\amin}{a_{\mathrm{min}}}
\newcommand{\cmin}{c_{\mathrm{min}}}

\newcommand{\amax}{a_{\mathrm{max}}}
\newcommand{\Id}{\mathbf{Id}}
\newcommand{\bxi}{\bvec{\xi}}

\newcommand{\bA}{\bmat{A}}
\newcommand{\bB}{\bmat{B}}
\newcommand{\bC}{\bmat{C}}
\newcommand{\bX}{\bmat{X}}
\newcommand{\ba}{\bvec{a}}
\newcommand{\bc}{\bvec{c}}
\newcommand{\bcau}{\bvec{c}^a_{u; \nb, \nq}}
\newcommand{\bcaueps}{\bvec{c}^a_{u; \Nbeps, \Nqeps}}
\newcommand{\bcf}{\bvec{c}_{f; \nb}}
\newcommand{\tbcf}{\widetilde{\bvec{c}}_{f; \nb}}

\newcommand{\bpsi}{\bm{\psi}}

\newcommand{\tp}{\widetilde{p}}

\newcommand{\bff}{\bvec{f}}
\newcommand{\ap}{\mathfrak{a}}
\newcommand{\tap}{\widetilde{\mathfrak{a}}}
\newcommand{\Cpoin}{C_{\mathrm{poi}}}
\newcommand{\Ccoer}{C_{\mathrm{coer}}}
\newcommand{\Ccont}{C_{\mathrm{cont}}}

\newcommand{\CinvA}{C_{\mathrm{inv}, A}}
\newcommand{\Ccu}{C_{c_u}}
\newcommand{\Cbasis}{C_\mathrm{b}}
\newcommand{\Cfe}{C_\mathrm{G}}
\newcommand{\bfe}{b_\mathrm{G}}
\newcommand{\Cpol}{C_\mathrm{pol}}
\newcommand{\CAinv}{C_{\mathrm{A}}}

\newcommand{\Phiinputone}{\Phi^{A, \alpha}_{\nb , \nq }}
\newcommand{\Phiinputadr}{\Phi^{(\argsADR), \alpha}_{\nb , \nq }}
\newcommand{\Phiinputel}{\Phi^{\ttA, \alpha}_{\nb , \nq }}
\newcommand{\Phiinputtwo}{\Phi^{\widetilde{A}, \mathrm{Id},\alpha}_{\nb , \nq }}
\newcommand{\PhiinvA}{\Phi^{\widetilde{A}}_{\mathrm{inv}; \epsinv, \nb }}
\newcommand{\PhiinvAepsu}{\Phi^{\widetilde{A}}_{\mathrm{inv}; {\epsu}/(\|f\|_{L^2(\Omega)}\nb^{3/2}\CAinv), \nb }}
\newcommand{\Phiinv}{\Phi^{1-\delta, \nb }_{\mathrm{inv}; \frac{\epsinv}{\alpha}}}
\newcommand{\Phicu}{\Phi^{c_u}_{\epsu, \nb }}
\newcommand{\Phicufinal}{\Phi^{c_u}_{\epsu, \Nbeps}}
\newcommand{\Phibasis}{\Phi^{\mathrm{b}}_{\epsbasis, \nb }}
\newcommand{\Phibasisfinal}{\Phi^{\mathrm{b}}_{\epsbasis, \Nbeps}}
\newcommand{\Phibranch}{\Phi^{\mathrm{br}}_{\epsilon}}
\newcommand{\tPhibranch}{\widetilde{\Phi}^{\mathrm{br}}_{\epsilon/3}}
\newcommand{\Phibranchnoeps}{\Phi^{\mathrm{branch}}}
\newcommand{\Phitrunk}{\Phi^{\mathrm{tr}}_{\epsilon}}
\newcommand{\tPhitrunk}{\widetilde{\Phi}^{\mathrm{tr}}_{\epsilon/3}}
\newcommand{\Phitrunknoeps}{\Phi^{\mathrm{trunk}}}

\newcommand{\Phitap}{\Phi^{\tap}_{\epsp, \np, \nq}}

\newcommand{\Phiaone}{\Phi^{a,\mathrm{coef}}_{\epsilon, \np}}
\newcommand{\Phiaoneepsp}{\Phi^{a,\mathrm{coef}}_{\epsp, \np}}
\newcommand{\epsinv}{\epsilon_{\mathrm{inv}}}
\newcommand{\epsu}{\epsilon_{u}}
\newcommand{\epsbasis}{\epsilon_{\mathrm{b}}}
\newcommand{\epsfe}{\epsilon_{\mathrm{G}}}

\newcommand{\Nbeps}{\nb(\epsilon)}
\newcommand{\Nqeps}{\nq(\epsilon)}

\newcommand{\btau}{\bm{\tau}}
\newcommand{\bsigma}{\bm{\sigma}}
\newcommand{\bepsilon}{\bm{\epsilon}}

\newcommand{\frab}{\mathfrak{b}}
\newcommand{\frav}{\mathfrak{v}}
\newcommand{\frad}{\mathfrak{d}}

\newcommand{\Vpsi}{\bmat{V}_{\nq, \np}}
\newcommand{\baNq}{\frab^a_{\nq}}
\newcommand{\AaNb}{\bmat{A}^a_{\nb}}

\newcommand{\AaNbNq}{\bmat{A}^a_{\nb,\nq}}
\newcommand{\AoneNbNq}{\bmat{A}^1_{\nb,\nq}}
\newcommand{\AtildeNbNq}{\widetilde{\bmat{A}}^a_{\nb,\nq}}
\newcommand{\uaNbNq}{u^a_{\nb, \nq}}
\newcommand{\uaNbNqeps}{u^a_{\Nbeps, \Nqeps}}
\newcommand{\ugenparNbNq}{u^{\genpar}_{\nb, \nq}}
\newcommand{\axq}{\bm{a}_{\nq}}

\newcommand{\conc}{{\raisebox{2pt}{\tiny\newmoon} \,}}
\newcommand{\sconc}{\odot}

\newcommand{\R}{\mathbb{R}}
\newcommand{\Rsym}{\mathbb{R}_{\mathrm{sym}}}
\newcommand{\Rsymdd}{\Rsym^{d\times d}}

\newcommand{\Q}{\mathbb{Q}}

\newcommand{\xenc}{\bx_{\mathrm{enc}}}
\newcommand{\nb}{n_{\mathrm{b}}}
\newcommand{\tnb}{\widetilde{n}_b}
\newcommand{\tX}{\widetilde{X}}
\newcommand{\np}{n_\mathrm{p}}
\newcommand{\dimp}{d_\mathrm{p}}
\newcommand{\nq}{n_\mathrm{q}}
\newcommand{\argsADR}{\bA, c}
\newcommand{\bADR}{\frab^{(\argsADR)}}
\newcommand{\SADR}{S^{\mathrm{rd}}}
\newcommand{\cEADR}{\mathcal{E}^{\mathrm{rd}}_{\xenc}}
\newcommand{\tcU}{{\mathcal{U}^{\mathrm{rd}}}}
\newcommand{\cDel}{{\mathcal{D}^{\mathrm{el}}}}
\newcommand{\cUel}{{\mathcal{U}^{\mathrm{el}}}}
\newcommand{\tcD}{{\mathcal{D}^{\mathrm{rd}}}}
\newcommand{\Sel}{S^{\mathrm{el}}}
\newcommand{\ttB}{\mathtt{B}}
\newcommand{\cEel}{\mathcal{E}^{\mathrm{el}}_{\xenc}}
\newcommand{\cRel}{\cR^{\mathrm{el}}}
\newcommand{\Cq}{C_{\mathrm{q}}}

\newcommand{\Cb}{C_{\mathrm{b}}}
\newcommand{\Ap}{A_{\mathrm{p}}}
\newcommand{\epsp}{\epsilon_{\mathrm{p}}}

\newcommand{\dpi}{d_{\mathrm{p},i}}
\newcommand{\ttA}{{\mathtt A}}
\newcommand{\eps}{\bepsilon}
\DeclareMathOperator{\relu}{ReLU}
\let\div\relax
\DeclareMathOperator{\div}{div}
\DeclareMathOperator{\grad}{grad}
\DeclareMathOperator{\card}{card}

\newcommand{\phioned}{\phi^{\mathrm{1d}}}
\newcommand{\chioned}{\chi^{\mathrm{1d}}}
\newcommand{\nuoned}{\nu^{\mathrm{1d}}}
\newcommand{\zetaoned}{\zeta^{\mathrm{1d}}}

\newcommand{\Hper}{H_\mathrm{per}}
\newcommand{\Wper}{W_\mathrm{per}}

\newcommand{\CLell}{C_{L^2}}

\headers{Exponential Convergence of Deep Operator Networks for Elliptic PDEs}{C.~Marcati and Ch.~Schwab}
\title{Exponential Convergence of Deep Operator Networks\\
      for Elliptic Partial Differential Equations}

    \author{Carlo Marcati\thanks{Seminar for Applied
        Mathematics, ETH Zurich, Switzerland. Current address: Dipartimento di Matematica ``Felice Casorati'',
    University of Pavia, Italy (\email{carlo.marcati@unipv.it})} \and Christoph Schwab\thanks{Seminar for Applied
  Mathematics, ETH Zurich, Switzerland
  (\email{christoph.schwab@sam.math.ethz.ch})}}

\begin{document}
\maketitle
\begin{abstract}
We prove exponential expression rates of deep operator networks (deep ONets) between infinite-di\-men\-sion\-al spaces that emulate the coefficient-to-solution map of linear, elliptic, second-order, divergence-form partial differential equations (PDEs). In particular, we consider problems set in $d$-dimensional periodic domains, $d=1, 2, \dots$, with analytic right-hand sides and coefficients. Our analysis covers diffusion-reaction problems, parametric diffusion equations, and certain elliptic systems 
such as linear isotropic elastostatics in heterogeneous materials.

In the constructive proofs, we leverage the exponential convergence of spectral collocation methods for boundary value problems whose solutions are analytic. In the present periodic and analytic setting, this follows from classical elliptic regularity. Within the branch and trunk ONet architecture of Chen and Chen \cite{ChenChen1993} and of Lu et al.~\cite{lu2020deeponet}, we construct deep ONets which emulate the coefficient-to-solution map to a desired accuracy in the $H^1$ norm, uniformly over the PDE coefficient set. The ONets have size $\mathcal{O}(\left|\log(\varepsilon)\right|^\kappa)$, where $\varepsilon>0$ is the approximation accuracy, for some $\kappa>0$ depending on the physical space dimension.
\end{abstract}
\begin{keywords}
 Operator networks, deep neural networks, exponential convergence, elliptic PDEs
\end{keywords}
\begin{AMS}
35J15, 65N15, 65N35, 68T07
\end{AMS}
\section{Introduction}
\label{sec:intro}
The application of numerical surrogates of solution operators to partial
  differential equations (PDEs) 
via algorithms of deep learning 
has recently received considerable attention.
See, e.g., \cite{cai2021physicsinformed,Li2020,lu2021comprehensive,Lanthaler2021} and the references there.
Also, expression and approximation rate bounds for such computable operator
surrogates have appeared in various settings, see, e.g. 
\cite{Kovachki2021,deng2021convergence,DeRyck2022}, and the references there.
In the present paper,
we construct deep operator network (ONet) emulations
of coefficient-to-solution maps for boundary value problems 
with linear, second order elliptic divergence-form operators.
In particular, 
we consider operator networks with rectified linear unit (ReLU) activation
and problems formulated in domains without boundary and with analytic right-hand
sides and coefficients.  
In this setting, we construct operator networks that approximate the
(nonlinear) coefficient-to-solution map with exponential accuracy 
in the corresponding function spaces. 
We bound---poly-logarithmically with respect to the energy norm of the
error---both
the size of the approximating network
and the number of sampling points where the coefficient is queried.
\subsection{Existing Results}
\label{sec:ExRes}
Deep neural networks (DNN)  have been 
employed increasingly in recent years in the numerical solution of 
differential equations in science and engineering. 
We refer to the survey \cite{cai2021physicsinformed} 
for uses and successes of DNN based numerical simulations in computational fluid
mechanics, and to \cite{Ruf2020} for their use in computational finance and 
computational option pricing. 
First uses of DNNs in numerical PDE 
solution in engineering and the sciences focused on leveraging DNNs 
for ``mesh-free'' solution approximation and representation (see, e.g., \cite{PiNNs,DeepRitz}),
with good success explained, to some extent, by \emph{approximation properties of DNNs in 
function spaces} (see, e.g.,~\cite{PETERSEN2018296,OPS20_2738,Opschoor2019,MOPS20_938,SchwabZech,EGPGSurvey})
in particular overcoming the curse of dimensionality (CoD) in high-dimensional approximation 
of PDE solution manifolds \cite{SchwabZech,Grohs2021},
of parametric PDEs and of PDEs on high-dimensional state spaces, 
as arising, e.g., in computational finance
(see \cite{Ruf2020,BeckEJentzen} and the references therein).  

  Reference \cite{Kutyniok2019} addressed
  the expression rate of ReLU NNs for the solution maps of parametric PDEs.
  The analysis in that paper proceeds through the DNN emulation of 
\emph{reduced bases} for the approximation of solutions of the PDEs.
The expression rate bounds obtained in \cite{Kutyniok2019} are subject to 
strong hypotheses on the DNN expressivity of reduced bases for the PDEs of interest.
The parameter sets (i.e., the domains of the solution operator)
considered in \cite{Kutyniok2019} are finite-dimensional;
this paper mostly concerns instead  the approximation of solution maps
between infinite dimensional spaces. We nonetheless show how expression rates
for finitely-parametric PDEs also follow from our main results, see 
Theorem \ref{th:parametric} and Remark \ref{remark:kutyniok}.

DNNs have been leveraged in \cite{deng2021convergence,lu2020deeponet,Lanthaler2021}
for the \emph{DNN emulation of data-to-solution operators for PDEs}.  
See also the review \cite{lu2021comprehensive}.
Here, previous investigations have focused on \emph{universality of NNs}
for operator approximation. 
The pioneering work \cite{ChenChen1993} 
established this for a certain type of NNs with a ``branch and trunk'' architecture,
which will also be used in the present work. 
While \cite{ChenChen1993} imposed strong compactness assumptions, 
more recently \cite{Lanthaler2021} extended these results
to certain settings without the compactness assumptions of \cite{ChenChen1993}.
In these papers, 
focus has been on emulating nonlinear maps, such as domain-to-solution, 
or coefficient-to-solution maps. 
For well-posed PDE problems, continuous dependence on the problem data implies that 
these maps are continuous, 
in the appropriate topologies on the data and the solution space.
We refer to \cite{kovachki2021neural,lu2020deeponet} and the references therein.
In these references, 
some theory explaining some of the numerically observed performances 
of NN emulation of nonlinear operators 
has been developed (see, e.g., \cite{Lanthaler2021,deng2021convergence,DeRyck2022}).
We also mention the analysis of \cite{Kovachki2021} for Fourier Neural
  Operators, a different kind of operator networks, introduced in \cite{Li2020}.

The \emph{convergence rate estimates} proved in these references indicate
that a) DNNs are capable of parsimonious numerical representations of the 
nonlinear, smooth data-to-solution maps for PDEs, and 
b) they are not prone to the CoD
in connection with the countable number of parameters due, e.g., to 
series representations of inputs in 
separable Banach spaces of possibly infinite dimension.
\subsection{Contributions}
\label{sec:Contr}
We construct DNN approximations of data-to-solution maps, 
so-called ``Operator Networks'' 
for linear, second order divergence-form elliptic PDEs
with non-ho\-mo\-ge\-neous coefficients and source terms.
We establish \emph{exponential expression rates} for these 
coefficient-to-solution operators for elliptic PDEs.

Our argument relies on 
analytic regularity for elliptic PDEs with analytic coefficients, 
on the a priori analysis of periodic spectral approximation of PDEs,
and on the error analysis of numerical quadrature in fully discrete spectral methods. 
We consider linear second order divergence-form elliptic boundary value 
problems with analytic, periodic coefficients, and (uniformly) analytic solutions, 
whose inputs and solutions
admit exponentially convergent spectral collocation approximations
from spaces of high-degree, periodically extendable polynomials.
Our results show that
neural networks can emulate accurately the (nonlinear) data-to-solution operator of Galerkin methods 
for the elliptic PDEs mentioned above with numerical integration. 
The operator networks we construct are composed of \emph{encoding}, \emph{approximation}, and
\emph{reconstruction} operators. 
In the encoding step, the input datum 
is queried on collocation points in the physical domain.
The approximation and reconstruction parts of the 
operator networks are composed of two neural networks,
one that approximates a polynomial basis, while the other maps point evaluations
of the diffusion coefficient to coefficients over the basis.

Our proof is constructive, based on ``NN emulation'' of (building blocks of) 
  a spectral method.
  Our focus is on providing an upper bound on the expression rate of the 
  Operator Network approximation of the coefficient-to-solution maps, 
  rather than to suggest a concrete algorithm to actually construct those networks.  
  Learning algorithms may be able to perform the numerical Operator Network 
  construction more efficiently.

For the sake of clarity of exposition, we develop this strategy for 
model, linear second order elliptic PDEs 
in divergence form, with inhomogeneous coefficients.
We then show, using the compositionality of NNs, how 
to include problems with parametric diffusion, typically arising in
computational uncertainty quantification.
Finally, we mention the minor modifications required for PDEs with reaction
coefficients and discuss in some detail 
ONet emulation of the coefficient-to-solution map for linear elasticity.

The exponential expression of data-to-solution maps proved in this
  manuscript is the first result of this kind for operator networks.
 It is based on exponential compression rates of encoders and decoders
     which are based on spectral approximations to leverage analyticity
     of input and output of the data-to-solution maps.
     Here, analyticity of the solution is a consequence of classical 
     elliptic regularity. 
     The strong compression of spectral encoders and decoders 
     facilitated by analyticity allows to compose
     ONets from approximate neural network inversion of small, but generally
     dense spectral Galerkin matrices.
     By construction, therefore, the rates we obtain are bounded by the
       exponential ones
       of the spectral method emulated in our proofs. Although the problems considered are
       rather simple and well suited for classical numerical schemes, this work
       is also committed to laying the foundation for analyzing more
       complex problems. 
     ONet constructions for finite regularity input and output pairs
     with considerably different input encoder and output decoder maps 
     differ substantially from the present construction.
     They are considered in \cite{HSZ}.
  The approximation of data-to-solution maps for similar elliptic PDEs
  has been analyzed with different techniques in \cite{deng2021convergence}
  under weaker regularity assumptions on the coefficients.
  These lines of argument yield lower expression rate bounds.
\subsection{Structure of this paper}
\label{eq:struct}
To fix a setting for developing our results,
we introduce in Section~\ref{sec:problem} a scalar, elliptic, isotropic diffusion equation. 
The
\emph{coefficient-to-solution operator} 
that will be the main target of approximation by neural networks is
also introduced in this section. 
Then, in Section \ref{sec:ONets}, we define
feed forward neural networks (with ReLU activation)
and operator networks with the branch and trunk architecture of \cite{ChenChen1993,lu2020deeponet},
that approximate maps between infinite dimensional spaces.
We conclude the section by defining some operations on networks that will
then be used for the approximation analysis.
In Section \ref{sec:regularity-pol}, we gather (classical) results on the
polynomial approximation of solutions to the elliptic problem.
The main results of this paper are then proved in Section \ref{sec:NN-appx}. In
Theorem \ref{th:deepONet}, we show the exponential convergence of the operator
net approximation of the coefficient-to-solution map for the elliptic isotropic
diffusion problem. We extend the analysis to parametric diffusion coefficients
in Theorem \ref{th:parametric}. 
In Section \ref{sec:generalizations}, we
extend our ONet approximation to further second order problems 
comprising reaction-diffusion with nonzero reaction coefficients and 
linear elastostatics.
Finally, in Section \ref{sec:conc}, we recapitulate our results and outline
perspective extensions of our theory to elliptic periodic homogenization and to
further forward and inverse problems.


\subsection{Notation}
\label{sec:notation}
%
We use standard notation and symbols:
$\N$ denotes the set of positive natural numbers 
$\N=\{1, 2, 3, \dots\}$ and $\N_0 = \{0\}\cup \N$.
We write vectors in lowercase boldface characters and matrices
in uppercase boldface characters.
We denote by $\| \bvec{a}\|_2$ the $\ell^2$-norm of a vector $\bvec{a}$, while for
any matrix $\bmat{A}$, we denote 
$\| \bmat{A} \|_2 = \sup_{\|\bvec{x}\|_2=1}\|\bmat{A}\bvec{x}\|_2$ 
its operator norm. By $\|\bmat{A} \|_0$ and $\| \bvec{x}\|_0$ we denote,
  respectively, the number of nonzero
  elements of a matrix $\bmat{A}$ and a vector $\bvec{x}$.
The spectrum of a matrix $\bA$ is written $\sigma(\bA)$.
For $n\in \N$, $\Id_n$ is the $n\times n$ identity matrix, while $\bzero_n$ is a
vector of zeros of size $n$.
When used between matrices, we denote by $\otimes$ the Kronecker product: given two matrices
$\bA\in \R^{m\times n}$ and $\bB\in \R^{p\times q}$, then $\bC = \bA\otimes\bB\in
\R^{mp\times nq}$, such that
\begin{equation*}
  \bC = \begin{bmatrix}
  \bA_{11} \bB & \cdots & \bA_{1n}\bB \\
             \vdots & \ddots &           \vdots \\
  \bA_{m1} \bB & \cdots & \bA_{mn} \bB
\end{bmatrix}.
\end{equation*}
Given two functions $v_1$, $v_2$, we instead denote by $v_1 \otimes v_2$ the
  function such that $(v_1 \otimes v_2)(x_1, x_2) = v_1(x_1)v_2(x_2)$.
We denote by $[\ba_1 | \dots | \ba_n ]$ the matrix with columns $\ba_1, \dots, \ba_n$.
We indicate by $\vec : \R^{m\times n} \to \R^{mn}$ and
$\matr:\R^{mn}\to\R^{m\times n}$ the vectorization and matricization operators,
such that $\matr(\vec(\bA)) = \bA$ for any matrix $\bA$. All results are
independent of the ordering of the vectorization operation; the dimensions of
the matricization operation will be clear from the context.
We denote by $\Rsym^{n\times n}$ the space of symmetric matrices of size
$n\times n$.
Given two matrices $\bA, \bB\in \R^{n\times n}$, we write $\bA:\bB = \sum_{i,j=1}^n\bA_{ij}\bB_{ij}$.

Let $d\in \N$. 
For $k\in \N_0$, $p\in [1, \infty]$, and a domain $D\subset\R^d$, 
we indicate by $W^{k, p}(D)$ the classical Sobolev spaces. 
$W^{k,p}_{\mathrm{loc}}(\R^d)$ indicates functions that are in $W^{k, p}(D)$
for any bounded subset $D$ of $\R^d$.
In the Hilbertian case $p=2$, we write $H^k(D)$; 
in addition, $L^p(D) = W^{0,p}(D)$ and $L^2(D) = H^0(D)$.
Given $Q=(0,1)^d$ and $\Omega = (\R/\Z)^d$, 
we denote
, for all $k\in\N_0$ and $p\in[1, \infty]$,
\begin{equation*}
W^{k,p}(\Omega) 
= 
\Wper^{k, p}(Q) \coloneqq 
\left\{ v|_Q : v\in W^{k,p}_{\mathrm{loc}}(\R^d)  \text{ and $v$ is $Q$-periodic} \right\},
\end{equation*}
i.e., the restriction to $Q$ of all functions in
  $W^{k,p}_{\mathrm{loc}}(\R^d)$ that are $Q$-periodic.
We denote by $(\cdot, \cdot)$ the $L^2$ scalar product in $Q$.

For $C>0$, define $\Hol(\Omega; C)$ as the 
set of functions $v$ that are real analytic in $\R^d$, 
periodic with period one in all coordinate directions, 
and such that 
  \begin{equation}
    \label{eq:analytic-uniform}
    \| v \|_{W^{k, \infty}(Q)} \leq C^{k+1}k!, \qquad \forall k\in \N_0.
  \end{equation}
Define furthermore the set of all real analytic functions in $\Omega$ as
$\Hol(\Omega) = \bigcup_{C>0} \Hol(\Omega; C)$. 
By the Arzelà-Ascoli theorem, 
the set $\Hol(\Omega; C)$ is compact in $L^\infty(\Omega)$. 
\section{Problem formulation}
\label{sec:problem}

We introduce the set of admissible 
diffusion coefficient data $\cD$:
for each coefficient $a\in \cD$
we assume ellipticity in the form that
there exist constants $\amin, \amax > 0$ such that 
\begin{equation}
\label{eq:abounds}
\forall \bx\in Q,\, \forall a\in \cD\: \quad 
  \amin \leq a(\bx) \leq \amax
\;.
\end{equation}
We also assume that all $a\in \cD$ are real analytic and $Q$-periodic, 
with uniform bounds on the radius of convergence of the Taylor series: 
there exists a constant $A_{\cD}>0$ such that
\begin{equation} \label{eq:abound}
  \cD \subset \Hol(\Omega; A_{\cD}).
\end{equation}
As it will be useful in the sequel, 
we define the Poincaré constant $\Cpoin>0$
such that
\begin{equation}
\label{eq:poincare}  
\| v - \frac{1}{|Q|}\int_\Omega v\|_{L^2(Q)} 
  \leq 
   \Cpoin \| \nabla v\|_{L^2(Q)}, \qquad \forall v\in H^1(\Omega)=\Hper^1(Q).
\end{equation}
The ellipticity hypotheses \eqref{eq:abounds} and the Poincaré inequality
\eqref{eq:poincare} imply that
for every $f\in L^2(\Omega)$ such that $\int_{Q} f = 0$,
and for each $a\in \cD$, 
the elliptic boundary value problem
\begin{equation}
  \label{eq:problem}
  -\nabla \cdot (a\nabla u^a)  = f \text{ in }\Omega
\end{equation}
admits a unique solution
\begin{equation*}
u^a\in X\coloneqq \left\{ v\in H^1(\Omega) : \int_{{Q}}v =0 \right\} \simeq \Hper^1(Q)/{\mathbb R}.
\end{equation*}
It satisfies the \emph{variational formulation}: 
given $a\in \cD$, find $u\in X$
such that 
\begin{equation}\label{eq:wProblem}
\frab^a(u,v) = (f,v) \quad \forall v\in X.
\end{equation}
Here, for a given $a\in \Hol(\Omega)$, 
the bilinear form 
$\frab^a(\cdot,\cdot): H^1(\Omega) \times H^1(\Omega) \to \mathbb{R}$
is given by 
\begin{equation*}
\frab^a(w,v) \coloneqq \int_{{Q}} \left(  a \nabla w \cdot \nabla v  \right).
\end{equation*}
In what follows, 
we assume the fixed source term $f\in \Hol(\Omega) \cap X$ to be given and,
for any $a\in \cD$, we denote by $u^a$
the unique solution of \eqref{eq:wProblem} for this choice of $f$.

We denote (still keeping the source term $f$ in \eqref{eq:wProblem} fixed) 
by $S$ the data-to-solution operator $a\mapsto u^a$ in \eqref{eq:wProblem}.
We let $\cU = S(\cD)$ the set of solutions of \eqref{eq:wProblem} corresponding 
to inputs from $\cD$. 
As shown in Lemma \ref{lemma:lip-sol} in Appendix \ref{sec:LipSolOp},
for fixed right source term $f$ in \eqref{eq:wProblem},
the data-to-solution map $S:L^\infty(\Omega) \to H^1(\Omega)$ is Lipschitz continuous.
Furthermore, standard elliptic regularity (see \cite{Morrey,Costabel2010}
and Lemma \ref{lemma:regularity} below) implies $S(\cD) \subset \Hol(\Omega)$.
\section{Neural operator networks}
\label{sec:ONets}
Our goal is to derive bounds for the approximation of the 
solution operator $S:\cD\to X\subseteq H^1(\Omega)$,
defined in Section \ref{sec:problem}, by an operator network.
To define operator networks, we recall the definition of
classical feed forward neural networks with ReLU activation
\begin{equation*}
\relu: \R \to \R, \, x \mapsto \max\{0, x\}.
\end{equation*}
\subsection{Feed forward neural network}
\label{sec:NN}

\begin{definition}[{\cite[Definition 2.1]{PETERSEN2018296}}] 
\label{def:NeuralNetworks}
Let $d, L\in \N$. 
A \emph{neural network $\Phi$ with input dimension $d$ and $L$ layers} 
is a sequence of matrix-vector tuples 
\[
    \Phi = \big((\bA_1,\bb_1),  (\bA_2,\bb_2),  \dots, (\bA_L, \bb_L)\big), 
\]
where $N_0 \coloneqq d$ and $N_1, \dots, N_{L} \in \N$, and 
where $\bA_\ell \in \R^{N_\ell\times N_{\ell-1}}$ and $\bb_\ell \in \R^{N_\ell}$
for $\ell =1,...,L$.

For a NN $\Phi$, 
we define the associated
\emph{realization of the NN $\Phi$} as 
\[
 \mathrm{R}(\Phi): \R^d \to \R^{N_L} , \, \bx\mapsto \bx_L \eqqcolon \mathrm{R}(\Phi)(\bx),
\]
where the output $\bx_L \in \R^{N_L}$ results from 
\begin{equation}
    \label{eq:NetworkScheme}
    \begin{split}
        \bx_0 &\coloneqq \bx, \\
        \bx_{\ell} &\coloneqq \relu(\bA_{\ell} \, \bx_{\ell-1} + \bb_\ell), \qquad \text{ for } \ell = 1, \dots, L-1,\\
        \bx_L &\coloneqq \bA_{L} \, \bx_{L-1} + \bb_{L}.
    \end{split}
\end{equation}
Here $\relu$ is understood to act component-wise on vector-valued inputs, 
i.e., for $\by = (y^1, \dots, y^m) \in \R^m$,  $\relu(\by) \coloneqq (\relu(y^1), \dots, \relu(y^m))$.
We call $N(\Phi) \coloneqq d + \sum_{j = 1}^L N_j$ the \emph{number of neurons of the NN} 
$\Phi$, $\depth(\Phi)\coloneqq L$ the \emph{number of layers} or 
\emph{depth}, $\size_j(\Phi)\coloneqq \| A_j\|_{0} + \| b_j \|_{0}$ 
the \emph{number of nonzero weights in the $j$-th layer}, 
and
$\size(\Phi) \coloneqq \sum_{j=1}^L \size_j(\Phi)$ the \emph{number of nonzero weights of $\Phi$}, 
also referred to as its \emph{size}. 
We refer to $N_L$ as the \emph{dimension of the output layer of $\Phi$}.
\end{definition}

\subsection{Operator networks}
\label{sec:subsec-ONet}
\begin{figure}
  \centering
  \includegraphics[width=.5\textwidth]{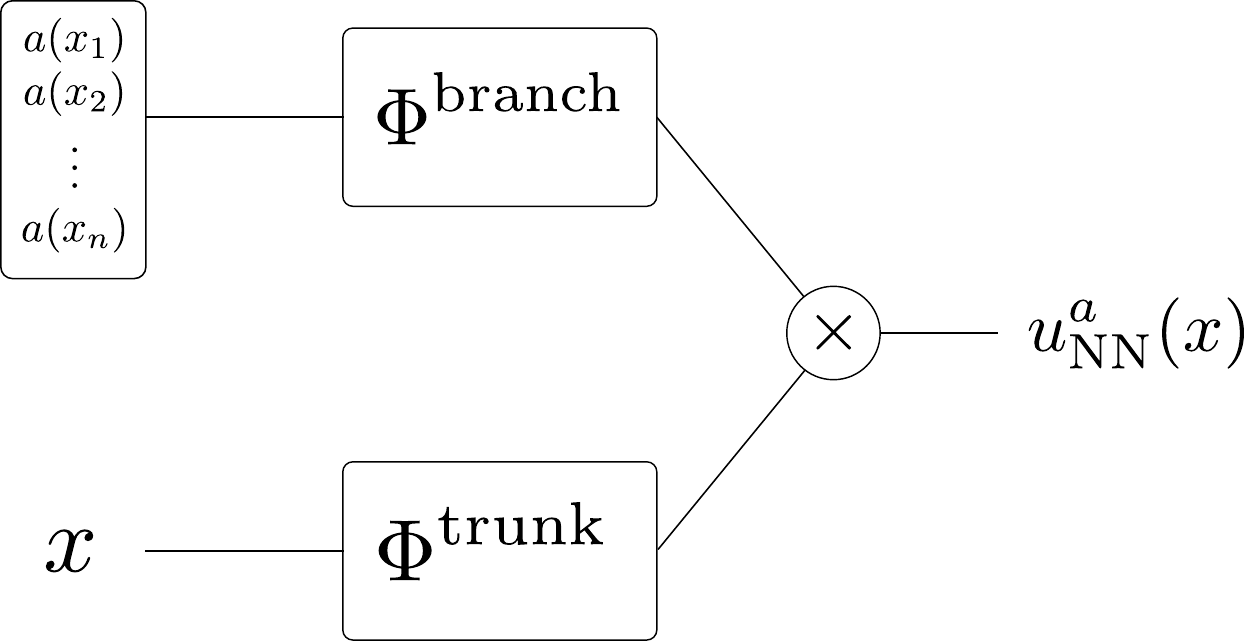}
  \caption{Structure of the branch and trunk network; $u^a_{\mathrm{NN}}(\bx) \coloneqq 
      \Realiz(\Phibranchnoeps) (\ba) \cdot \Realiz(\Phitrunknoeps)(\bx)$.}
  \label{fig:ONet}
\end{figure}

The operator network approximating the solution operator $S$ can be seen as the
composition $\mathcal{R}\circ \mathcal{A} \circ \mathcal{E}$ of three 
mappings: 
\begin{itemize}
\item Encoding $\mathcal{E}: \cD \to \R^n$, for $n\in \N$,
\item Approximation $\mathcal{A}: \R^n\to \R^m$, for $m\in\N$,
\item Reconstruction $\mathcal{R}: \R^m \to H^1({Q})$,
\end{itemize}
see the diagram in Figure \ref{fig:ONet-scheme}. 
We refer the reader to \cite{Lanthaler2021,kovachki2021neural} 
for a broader view on and thorough discussion of operator networks between
infinite dimensional spaces.
\begin{figure}
  \centering
  \begin{tikzpicture}
  \matrix (m) [matrix of math nodes,row sep=5em,column sep=8em,minimum width=2em]
  {\cD  &  H^1(\Omega)\\
     \mathbb{R}^n & \mathbb{R}^m\\};
  \path[-stealth]
    (m-1-1) edge node [left, align=center] {Encoding\\$a \mapsto \ba = [a(\bx_1), \dots, a(\bx_n)]$\\} (m-2-1)
            edge [dashed] node [above] {$S$} (m-1-2)
    (m-2-1.east|-m-2-2) edge node [below, align=center] {Approximation \\$\ba \mapsto \bc= [c_1, \dots , c_m]$} (m-2-2)
    (m-2-2) edge node [right, align=center] {Reconstruction \\ $\bc\mapsto u^a_{\mathrm{NN}} = \sum_i c_i \psi_i$} (m-1-2);
\end{tikzpicture}
  \caption{Diagram of operator network between infinite dimensional spaces}
  \label{fig:ONet-scheme}
\end{figure}
In our analysis, the encoding step will map functions $a\in\cD$ 
to the vector $\ba\in \mathbb{R}^n$ of their point evaluations, i.e.
\begin{equation*}
\ba = \cE_{\{\bx_1, \dots, \bx_n\}}(a)\coloneqq \left[  a(\bx_1), \dots, a(\bx_n) \right]^\top, 
\end{equation*}
for suitable collection of points $\bx_1, \dots, \bx_n \in \overline{Q}$.
The approximation and reconstruction steps involve feed-forward neural
  networks that we will, respectively, denote $\Phibranchnoeps$ and
  $\Phitrunknoeps$. Specifically, the \emph{approximate solution operator}
$\mathcal{A}$ is realized as 
\begin{equation*}
  \mathcal{A}_{\Phibranchnoeps}(\ba) = \Realiz(\Phibranchnoeps)(\ba).
\end{equation*}
For the reconstruction step $\cR$,
for all $\bc\in \mathbb{R}^m$ and $\bx\in \overline{Q}$, we define
\begin{equation*}
  \cR_{\Phitrunknoeps}(\bc) (\bx) =  \left( \Realiz(\Phitrunknoeps) (\bx)  \right)\cdot \bc. 
\end{equation*}
This constructs the operator network mapping from $\cD$ to $H^1({Q})$, defined by
\begin{equation*}
  \cR_{\Phitrunknoeps}\circ \cA_{\Phibranchnoeps} \circ \cE_{\{\bx_1, \dots, \bx_n\}}: 
  a\mapsto u^a_{\mathrm{NN}}(\cdot) \coloneqq
  (\Realiz(\Phitrunknoeps)(\cdot) )\cdot\Realiz(\Phibranchnoeps)(\mathcal{E}_{\bx_1, \dots, \bx_n}(a))
\end{equation*}
see Figure \ref{fig:ONet}. 
For the precise definition of the branch and trunk networks used to
  approximate the solution operator of \eqref{eq:problem} we refer the
  reader to Sections \ref{sec:branch} and \ref{sec:trunk}.

We aim for operator networks that approximate, for all $a\in \cD$,
solutions $u^a$ of \eqref{eq:problem} in the $H^1({Q})$-norm, 
uniformly over the input space $\cD$, i.e., such that
\begin{equation*}
  \sup_{a\in \cD}\| u^a - u^a_{\mathrm{NN}} \|_{H^1({Q})} \leq \epsilon .
\end{equation*}
The main result of this paper consists in proofs for upper bounds on 
$n$, $m$, and on the sizes of $\Phitrunknoeps$ and $\Phibranchnoeps$ 
as functions of the error $\epsilon$.
\subsection{Operations on neural networks}
\label{S:OpNN}
We introduce and recall some operations on neural networks that
will be necessary for the construction of the branch and trunk networks.
\subsubsection {Concatenation and sparse concatenation}
\begin{definition}[NN concatenation, {{\cite[Definition 2.2]{PETERSEN2018296}}}]
Let $L_1, L_2 \in \N$ and let
\[\Phi^1 = ((\bA_1^1,\bb_1^1), \dots, (\bA_{L_1}^1,\bb_{L_1}^1)), \qquad
\Phi^2 = ((\bA_1^2,\bb_1^2), \dots, (\bA_{L_2}^2,\bb_{L_2}^2))
\]
be two neural networks such that the input layer of $\Phi^1$ has the same
dimension as the output layer of $\Phi^2$.
Then, \emph{$\Phi^1 \conc \Phi^2$} denotes the following $L_1+L_2-1$ layer
network:
\[
  \Phi^1 \conc \Phi^2
  \coloneqq  ((\bA_1^2,\bb_1^2),
      \dots,
      (\bA_{L_2-1}^2,\bb_{L_2-1}^2),
      (\bA_{1}^1 \bA_{L_2}^2, \bA_{1}^1 \bb^2_{L_2} + \bb_1^1),
      ({\bA}_{2}^1, \bb_{2}^1),
      \dots,
      (\bA_{L_1}^1, \bb_{{L_1}}^1)
     ).
\]
We call $\Phi^1 \conc \Phi^2$ the \emph{concatenation of $\Phi^1$ and $\Phi^2$}.
\end{definition}

\begin{proposition}[NN sparse concatenation, {{\cite[Remark 2.6]{PETERSEN2018296}}}]
\label{prop:conc}
Let $L_1, L_2 \in \N$, and let 
$\Phi^1, \Phi^2$ 
be two NNs of respective depths $L_1$ and $L_2$ such that $N^1_0 = N^2_{L_2}\eqqcolon d$, i.e.,
the input layer of $\Phi^1$ has the same dimension as the output layer of $\Phi^2$. 

Then, there exists a NN $\Phi^1 \sconc \Phi^2$, called 
the \emph{sparse concatenation of $\Phi^1$ and $\Phi^2$}, 
such that $\Phi^1 \sconc \Phi^2$ has $L_1+L_2$ layers,   
$\Realiz(\Phi^1 \sconc \Phi^2) = \Realiz(\Phi^1) \circ \Realiz(\Phi^2)$ 
and 
$\size\left(\Phi^1 \sconc \Phi^2\right) \leq 2\size\left(\Phi^1\right)  + 2\size\left(\Phi^2\right)$.
\end{proposition}
\subsubsection{Emulation of matrix inversion}
\label{sec:MatInvReLU}
Dense matrix inversion can be approximated by suitable ReLU NNs.
We recall the following result from \cite{Kutyniok2019} where,
for $Z\in \mathbb{R}_+$ and $N\in\N$, 
\begin{equation*}
  K^Z_N \coloneqq \{\vec(\bA) : \bA\in \mathbb{R}^{N\times N} , \|\bA\|_2 \leq Z\}.
\end{equation*}
\begin{theorem}{{\cite[Theorem 3.8]{Kutyniok2019}}} \label{th:inversion}
 For $\epsilon,\delta\in (0,1)$ define
\[
m(\epsilon,\delta)\coloneqq \left\lceil \frac{\log\left( 0.5 \epsilon\delta \right)}{\log(1-\delta)}\right\rceil.
\]
There exists a universal constant $C_{\mathrm{inv}} > 0$ 
such that 
for every $N\in \N,$ $\epsilon\in (0,1/4)$ and every $\delta\in (0,1)$ 
there exists  a NN $\Phi^{1-\delta, N}_{\mathrm{inv},\epsilon}$ 
with $N^2$-dimensional input, $N^2$-dimensional output 
and the following properties:
\begin{enumerate}
\item $\depth\left(\Phi^{1-\delta, N}_{\mathrm{inv};\epsilon}\right)\leq C_{\mathrm{inv}} \left(1+\log\left(m(\epsilon,\delta) \right)  \right)\cdot  \left(\log\left(1/ \epsilon\right)+\log\left(m(\epsilon,\delta) \right)+\log(N)\right)$,
\item $\size\left(\Phi^{1-\delta, N}_{\mathrm{inv};\epsilon}\right)\leq C_{\mathrm{inv}} m(\epsilon,\delta) \left(1+\log^2(m(\epsilon,\delta))  \right)N^3\cdot  \left(\log\left(1\epsilon\right)+ \log\left(m(\epsilon,\delta)\right) +\log(N) \right)$,
    
    \item $\sup_{\vec(\bA)\in K^{1-\delta}_{N}} \left\|\left(\Id_{{N}}-{\bA} \right)^{-1} - \matr\left( \Realiz\left(\Phi^{1-\delta,N}_{\mathrm{inv};\epsilon}\right)(\vec(\bA))\right)   \right\|_2 \leq \epsilon,$
    \item for any $\vec(\bA)\in K^{1-\delta}_N$ we have 
    \begin{align*}
        \left\| \matr\left( \Realiz\left(\Phi^{1-\delta,N}_{\mathrm{inv};\epsilon}\right)(\vec(\bA))\right)   \right\|_2\leq \epsilon + \left\|\left({\Id}_{N}-\bA\right)^{-1} \right\|_2 \leq \epsilon+ \frac{1}{1-\|\bA\|_2} \leq \epsilon + \frac{1}{\delta}.
    \end{align*}
\end{enumerate}
\end{theorem}
\begin{remark}
  The bounds on the depth and size of the network of
    Theorem \ref{th:inversion} are slightly modified 
    compared to those in \cite{Kutyniok2019}, since some instances of
    $\log(m(\epsilon, \delta))$ have been replaced by $1+\log(m(\epsilon,
    \delta))$. Indeed, for all $\delta \in [2/(2+\epsilon), 1)$, with fixed $\epsilon>0$,
    $m(\epsilon,\delta) = 1$.
    In this case, the unmodified estimates would give a
    degenerate bound on the depth and size of the network. 
    This modification is mathematically inconsequential, 
    the relevant case for the approximation estimates being $\epsilon\downarrow 0$.

\end{remark}
\section{Regularity and polynomial approximation}
\label{sec:regularity-pol}
We shall exploit the classical fact that the analyticity of the coefficient $a$ 
    and of the source term $f$ in $\Omega$ combined with periodicity
    implies 
    analyticity of the solution $u^a$ of \eqref{eq:problem}.
    This, in turn, will imply exponential convergence of tensor product polynomial 
    (spectral) approximations of $a$ and $u^a$, 
    which will be the basis of the NN approximation 
    developed in Section~\ref{sec:NN-appx} ahead.
\subsection{Regularity}
\label{sec:regularity}
The following result follows from \cite[Remark 1.6.5 and Theorem 1.7.1]{Costabel2010}.
\begin{lemma}
  \label{lemma:regularity}
  There exists $A_{\cU}> 0$ such that $S(\cD)\subset \Hol(\Omega; A_{\cU})$.
\end{lemma}
\begin{proof}
  From \cite[Theorem 1.7.1]{Costabel2010}, it follows that
  $S(\cD)\subset \Hol(\Omega)$ and, for each $u\in S(\cD)$, there exists $A_u>0$
  such that
  \begin{equation*}
    \frac{1}{k!}| u |_{H^{k}({Q})} \leq A_u^{k+1}\left( \sum_{j=0}^{k-2} \frac{1}{j!} |f|_{H^{j}({Q})} +
   \| u\|_{H^1({Q})})\right), \qquad \forall k\in \N_0.
  \end{equation*}
  Furthermore, from \cite[Remark 1.6.5]{Costabel2010}, inspecting the proof
  of \cite[Theorem 1.7.1]{Costabel2010}, and from \eqref{eq:abound}, the proof is completed 
  since it follows that 
  \begin{equation*}
    A_{\cU} \coloneqq \sup_{u\in S(\cD)} A_u  < \infty.
  \end{equation*}
\end{proof}
\subsection{Polynomial basis and quadrature}
\label{sec:pol-quad}
Consider the univariate Legendre polynomials $L_0$, $L_1$, $\dots$ such that $L_i\in\Q_i((0,1))$,
    normalized with $L_i(1)=1$. 
Define then, for all $i \in \N $,
    \begin{equation*}
      \phioned_0 = L_0, \qquad
      \phioned_{2i-1} = L_{2i}, \qquad \phioned_{2i} = L_{2i+1} - L_1.
    \end{equation*}
    These functions satisfy, for all $i\in\N_0$, $\phioned_i(0) = \phioned_i(1)$.
    It follows that, for all $p\in\N$,
    \begin{equation*}
      \spn(\phi_1, \dots, \phi_{p}) = \Q_{p+1}((0,1)) \cap \left\{v \in
      H^1(\R/\Z): \int_{{(0,1)}}v = 0\right\}.
    \end{equation*}
    We can then introduce, for all integer $p\geq 2$,
    \begin{equation}
      \label{eq:phibasis}
      \phi_{i_1 + p i_2 + \dots + p^{d-1}i_d} = \phioned_{i_1} \otimes \phioned_{i_1} \otimes \dots \otimes \phioned_{i_d}, \qquad (i_1, \dots, i_d)\in \{0, \dots, p-1\}^d.
    \end{equation}
    Then, as shown in
      Lemma \ref{lemma:basis} in the appendix, for all integer $p\geq 2$, denoting $\nb = p^d -1$, 
    \begin{equation}
      \label{eq:Xnb}
      X_{\nb} \coloneqq \spn(\{\phi_1, \dots, \phi_{\nb}\}) = \left\{ v\in \Q_{p}(Q) : \int_{Q}v = 0 \text{ and } v\in H^1(\Omega) \right\} =  \Q_p(Q)\cap X.
    \end{equation}
    The restriction to polynomials of degree $p\geq 2$ is without loss of generality,
    as the periodicity and vanishing average constraints imply $\Q_1(Q) \cap X = \{0\}$.

For a quadrature order parameter  $q\geq 2$, denoting $\nq =q^d$, 
we consider the
\emph{Gauss-Lobatto quadrature rule} 
with weights $\{\wkq\}_{k=1}^{\nq }$ 
and points
$ \{\xkq\}_{k=1}^{\nq }\subset \overline{Q}$ such that
\begin{equation*}
  \int_{{Q}} g = \sum_{k=1}^{\nq } \wkq g(\xkq), 
  \qquad \forall g\in \mathbb{Q}_{2q-3}(Q).
\end{equation*}
There exist constants $\cqone, \cqtwo>0$ such that
\begin{equation}
  \label{eq:quad-equiv}
  \cqone \| v \|^2_{L^2({Q})} \leq \sum_{k=1}^{(p+1)^d} \wkppone (v(\xkppone))^2 
  \leq \cqtwo \| v\|^2_{L^2({Q})},\quad \forall v\in \Q_{p}(Q), \, \forall p\in \N,
\end{equation}
  see \cite[Equation (6.4.52)]{CHQZ1}. 
  We remark that the constants $\cqone$ and $\cqtwo$ are independent of $p$, 
  but depend in general exponentially on the dimension $d$.
We may assume, without loss of generality, that $\cqone \leq 1$ and $\cqtwo \geq 1$.
We introduce furthermore the bilinear form with quadrature $\baNq$ 
\begin{equation*}
\baNq(u, v) = \sum_{k=1}^{\nq } \wkq a(\xkq) \nabla u(\xkq) \cdot \nabla v (\xkq), 
\qquad \forall u, v\in C^1(\overline{Q}).
\end{equation*}
Eventually, here $u,v$ shall be tensor product polynomials in $Q$.

For each $a\in \cD$, we introduce the symmetric matrices 
\begin{equation*}
  \left[\AaNb   \right]_{ij} = \frab^a(\phi_j, \phi_i), 
  \qquad
  \left[\AaNbNq   \right]_{ij} = \baNq(\phi_j, \phi_i),
  \qquad
  (i, j)\in\brange{\nb }^2.
\end{equation*}
Let $\AoneNbNq$ be the matrix obtained with $a \equiv 1$ in ${Q}$.
Let $q\geq p+1$:
  for all nonzero $\bx\in\R^{\nb}$, 
  there exists $v\in X_{\nb}\backslash\{ 0 \}$
  such that, for all $a\in \cD$,
  \begin{equation}
    \label{eq:Aposdef}
    \bx^\top \AaNbNq \bx = \baNq(v, v) > 0, 
  \end{equation}
  due to the equivalence of norms \eqref{eq:quad-equiv} and to the 
  Poincar\'{e} inequality \eqref{eq:poincare}.
  Hence, the matrices $\AaNbNq$ and $\AoneNbNq$ are invertible.
Denote then
\begin{equation*}
  \AtildeNbNq = \rev{( \AoneNbNq )^{-1/2} \AaNbNq( \AoneNbNq )^{-1/2}}.
\end{equation*}
We also introduce the right-hand side vector $\bcf\in \mathbb{R}^{\nb }$ 
such that
\begin{equation}
  \label{eq:bcf-def}
  \left[ \bcf \right]_i = \int_{{Q}} f \phi_i, \qquad i\in\brange{\nb }\;.
\end{equation}
The Cauchy-Schwarz inequality and
  \begin{equation*}
    \|\phi_i\|^2_{L^2(Q)}
    \leq 
    \begin{cases}
      1/(2i+3) &\text{if $i$ is odd}\\
      1/(2i+3) + 1/3 &\text{if $i$ is even},
    \end{cases}
     \qquad \forall i\in\N,
  \end{equation*}
  hence $\| \phi_i \|_{L^2(Q)}\leq 1$, imply that
  \begin{equation}
    \label{eq:bound-ctof}
    \| \bcf \|_{2}^2 \leq \|f\|_{L^2(Q)}^2 \sum_{i=1}^{\nb} \|\phi_i\|_{L^2(Q)}^2 \leq \nb \|f\|_{L^2(Q)}^2.
  \end{equation}
Finally, let
  \begin{equation*}
  \tbcf = (\AoneNbNq)^{-\rev{1/2}}\bcf.
  \end{equation*}%
  and
\begin{equation}
  \label{eq:cudef}
  \bcau \coloneqq (\AaNbNq)^{-1}\bcf \rev{=( \AoneNbNq )^{-1/2}(\AtildeNbNq)^{-1}\tbcf}.
\end{equation}
Here and in the sequel, for all $q \in \N$, with $\nq=q^d$, we will denote
$\axq \in \mathbb{R}^{\nq }$ the vector 
with entries
\begin{equation*}
  \left[\axq  \right]_i = a(\xkq), \qquad \forall k \in \brange{\nq }.
\end{equation*}

The following two statements concern the norms of the matrices introduced, and
will be necessary for later estimates.
\begin{lemma}
  \label{lemma:contcoer}
  There exist $\Ccoer, \Ccont>0$ such that for all $a\in \cD$,
  all $p\in\N$, and all integer $q \geq p+1$,
  \begin{equation}
    \label{eq:sigmaA}
    \sigma( \AtildeNbNq ) \subset [\Ccoer,\Ccont], \qquad  \sigma(( \AtildeNbNq )^{-1})\subset \left[{1}/{\Ccont}, {1}/{\Ccoer}  \right],
  \end{equation}
  with $\nb = p^d -1$ and $\nq = q^d$.
\end{lemma}
\begin{proof}
  For all $\bx\in \R^{\nb}$,
  \begin{equation*}
    \amin \bx^\top\AoneNbNq \bx\leq \bx^\top\AaNbNq \bx\leq \amax \bx^\top\AoneNbNq \bx.
  \end{equation*}
  Since the matrices $\AoneNbNq$ and $\AaNbNq$ are symmetric and positive
    definite, see \eqref{eq:Aposdef}, this implies, by Lemma
    \ref{lemma:eigprob} in the Appendix,
  \begin{equation*}
    \sigma(\rev{\AtildeNbNq }) \subset [\amin , \amax].
  \end{equation*}
  The assertion follows with $\Ccoer = \amin$ and $\Ccont = \amax$.
\end{proof}
We assume, for ease of notation, that $\Ccoer\leq 1$ and $\Ccont \geq 1$.
\begin{lemma}
  \label{lemma:norminvA}
  There exists a constant $\CAinv>0$ such that, for all $p\in\N$, and for all integer $q\geq p+1$,
  \begin{equation*}
    \| (\AoneNbNq )^{-1}\|_2\leq \CAinv\nb,
  \end{equation*}
  with $\nb = p^d -1$ and $\nq = q^d$.
\end{lemma}
\begin{proof}
  From \eqref{eq:quad-equiv}, \eqref{eq:poincare},
  the symmetry of the bilinear form, and Lemma \ref{lemma:L2l2}, it follows that
 \begin{align*}
   \lambda_{\min{}}(\AoneNbNq) &\coloneqq \min_{\lambda\in \sigma(\AoneNbNq)}\lambda 
   = \inf_{\bx\in \R^{\nb}} \frac{\bx^\top \AoneNbNq \bx}{\|\bx\|_2^2}
   \overset{\eqref{eq:L2-l2}}{\geq}
   \inf_{v\in X_{\nb}}\frac{\frab^1_{\nq}(v, v)}{\CLell^2\nb\|v\|_{L^2({Q})}^2} 
   \\
   \overset{\eqref{eq:quad-equiv}}&{\geq} \frac{\cqone^2}{\CLell^2\nb} \inf_{v\in X_{\nb}} 
   \frac{\|\nabla v\|_{L^2({Q})}^2}{\|v\|_{L^2({Q})}^2}
   \overset{\eqref{eq:poincare}}{\geq} \frac{\cqone^2}{\Cpoin^2\CLell^2\nb}.
 \end{align*}
 This concludes the proof, since
 $ \| (\AoneNbNq)^{-1} \|_2  = 1 / \lambda_{\min{}}(\AoneNbNq)$.
\end{proof}
\section{NN approximation}
\label{sec:NN-appx}
We detail the structure of the branch and trunk networks (see Fig.~\ref{fig:ONet})
and state and prove their convergence rate bounds.
%
\subsection{Branch network}
\label{sec:branch}
\subsubsection{Input layer}
\label{sec:InpLayer}
For all $k\in \brange{\nq }$, 
let $\hbD_{\nb }(\xkq)$ denote the matrix with entries
\begin{equation*}
  \left[\hbD_{\nb }(\xkq)  \right]_{ij} 
   = \wkq \nabla\phi_i(\xkq)\cdot \nabla\phi_j(\xkq), \qquad (i,j)\in \brange{\nb }^2.
\end{equation*}
The following statement follows from this definition. 
\begin{lemma}
  \label{lemma:Phiinputone}
For all $\alpha \in \mathbb{R}$, the one-layer NN 
\begin{equation*}
  \Phiinputone
  \coloneqq \left(
    \left(-\alpha \left[\vec(\hbD_{\nb }(\xq_1)) | \dots |\vec(\hbD_{\nb }(\xq_{\nq }))\right], \bzero_{\nb^2} \right)
  \right)
\end{equation*}
satisfies $\size(\Phiinputone) \leq \nb ^2\nq  $ and is such that 
\begin{equation}
  \label{eq:Phiinputone-thesis}
  \matr\left(\Realiz(\Phiinputone) (\axq)   \right) 
   = 
  -\alpha \AaNbNq.
\end{equation}
\end{lemma}
\begin{proof}
  We have
  \begin{equation*}
    \left[\Realiz\left( \Phiinputone \right) (\axq)   \right]_i= -\alpha\sum_{k=1}^{\nq }\left[\vec(\hbD(\xkq))\right]_i a(\xkq) ,
  \end{equation*}
  hence the equality after matricization. 
  The size bound follows from the fact that 
  \begin{equation*}
    \| \hbD_{\nb }(x^{(q)}_k) \|_0 \leq \nb^{2} , \;\; k\in \brange{\nq }\;.
  \end{equation*}
\end{proof}
\begin{lemma}
 \label{lemma:Phiinputtwo} 
For all $\alpha \in \mathbb{R}$, the \rev{three}-layer NN
\begin{equation}
  \label{eq:Phiinputtwo-def}
  \Phiinputtwo \coloneqq
  \left( \left( \Id_{\nb}\otimes (\AoneNbNq)^{-\rev{1/2}} , \vec(\Id_{\nb})\right) \right) \sconc
 \rev{ \left( \left( (\AoneNbNq)^{-1/2}\otimes \Id_{\nb} , \vec(\bzero_{\nb})\right) \right)\sconc}
  \Phiinputone
\end{equation}
is such that
\begin{equation*}
  \matr\left(\Realiz(\Phiinputtwo) (\axq)   \right) = 
  -\alpha \AtildeNbNq + \Id_{\nb}, 
\end{equation*}
and 
$\size(\Phiinputtwo) \leq \rev{6}\nb( 1 + \nb \nq  + \nb^2)$.
\end{lemma}
\begin{proof}
  For all $\rev{n}\in \N$ and all $\bA\in \R^{\rev{n\times n}}$ and
  $\bX\in\R^{\rev{n\times n}}$, $ \vec(\bA\bX) = (\Id_{\rev{n}}\otimes \bA) \vec(\bX)$
  \rev{and $\vec(\bX\bA^\top) = (\bA\otimes\Id_n)\vec(\bX)$}.
  Hence,
  \begin{align*}
    \Realiz(\Phiinputtwo)
    \overset{\eqref{eq:Phiinputtwo-def}}&{=}
   (\Id_{\nb}\otimes (\AoneNbNq)^{-\rev{1/2}})\rev{((\AoneNbNq)^{-1/2}\otimes \Id_{\nb})}\Realiz(\Phiinputone) + \vec(\Id_{\nb})\\
   \overset{\eqref{eq:Phiinputone-thesis}}&{=} -\alpha (\Id_{\nb}\otimes (\AoneNbNq)^{-\rev{1/2}})\rev{((\AoneNbNq)^{-1/2}\otimes \Id_{\nb})}\vec(\AaNbNq)+ \vec(\Id_{\nb})\\
   & = -\alpha \vec((\AoneNbNq)^{-\rev{1/2}}\AaNbNq\rev{(\AoneNbNq)^{-\rev{1/2}}})+ \vec(\Id_{\nb}).
  \end{align*}
  Finally, by Proposition \ref{prop:conc},
  \begin{align*}
    \size(\Phiinputtwo)
    &\leq \rev{4}\size(\Phiinputone) + \rev{6} \| \Id_{\nb}\otimes (\AoneNbNq)^{-\rev{1/2}} \|_0 + \rev{2}\|\Id_{\nb}\|_0
     \leq \rev{6}\size(\Phiinputone) + \rev{6} (\nb^3 + \nb).
  \end{align*}
\end{proof}
The following statement is a consequence of 
Lemmas \ref{lemma:contcoer} and \ref{lemma:Phiinputtwo}.
\begin{lemma}
  \label{lemma:inputnorm}
  Let $\Ccont, \Ccoer$ be the constants introduced in Lemma \ref{lemma:contcoer}.
 For all $a\in \cD$, 
 for all $p\in\N$, 
 for all integer $q\geq p+1$, and 
 for all $\alpha \in (0, 1/\Ccont )$,
  \begin{equation*}
    \| \matr\left(\Realiz (\Phiinputtwo) (\axq)\right) \|_2 \leq 1-\alpha\Ccoer \eqqcolon 1-\delta,
  \end{equation*}
  with $\nb =p^d-1$ and $\nq=q^d$.
\end{lemma}
\begin{proof}
  By Lemma \ref{lemma:Phiinputtwo}, $\matr\left(\Realiz (\Phiinputtwo)
    (\axq)\right) =  \Id_{\nb}-\alpha \AtildeNbNq  $.
Due to Lemma \ref{lemma:contcoer} and since $\alpha\leq 1/\Ccont$, this matrix
is symmetric positive definite and
\begin{equation*}
  \| \Id_{\nb}-\alpha \AtildeNbNq  \|_2  = \sup_{\substack{\bx \in \R^{\nb} \\ \|\bx\|_2=1}} \bx^\top (\Id_{\nb}-\alpha \AtildeNbNq )\bx = 1 - \alpha \inf_{\substack{\bx \in \R^{\nb} \\ \|\bx\|_2=1}}\bx^\top \AtildeNbNq \bx \leq 1-\alpha\Ccoer,
\end{equation*}
where the last inequality follows from Lemma \ref{lemma:contcoer}.
\end{proof}
Thanks to Theorem \ref{th:inversion} we can construct the network that approximates the inversion of the
``preconditioned'' Galerkin-Numerical Integration matrix $\AtildeNbNq$ 
 (more precisely, the network that emulates the map $\axq \mapsto ( \AtildeNbNq )^{-1}$).
\begin{proposition}
  \label{prop:invA}
  Let $\Ccoer, \Ccont$ be defined as in Lemma \ref{lemma:contcoer}. There
  exists a constant $\CinvA>0$ such that for all
  $\nb \in\N$ and for all $\epsinv \in (0, 1)$, 
  writing $\alpha = 1/(\Ccoer + \Ccont )$, 
  $\delta =\alpha\Ccoer$, $\nq  = \nb+1$, 
  and denoting
  \begin{equation}
\label{eq:PhiinvA-def}
    \PhiinvA \coloneqq \left( (\alpha \Id_{\nb }, \bzero_{\nb }) \right) \conc \Phiinv \sconc \Phiinputtwo,
  \end{equation}
  we have
  \begin{equation*}
    \sup_{a\in \cD}\|(\rev{\AtildeNbNq})^{-1} - \matr(\Realiz(\PhiinvA))(\axq) \|_2 \leq \epsinv, 
  \end{equation*}
  and  
  \begin{align*}
    \depth(\PhiinvA)&
    \begin{multlined}[t][.8\textwidth]
       \leq \CinvA \big[1+\log(1+\left|\log\epsinv\right|) + \log(\nb ) \big]\\
      \times \big[ 1+\left|\log\epsinv\right|  + \log(\nb ) +  \log(1+\left|\log\epsinv\right|) \big]
    \end{multlined}
    \\
    \size(\PhiinvA)&
                     \begin{multlined}[t][.8\textwidth]
                     \leq\CinvA  \nb ^{3}
\big[1+\left|\log\epsinv\right|\big] \big[ 1+ \log(1+\left|\log\epsinv\right|) + \log(\nb )\big]^2 
\\ \times\big[1+ \left|\log\epsinv\right| + \log(\nb ) +  \log(1+\left|\log\epsinv\right|) \big]\;.
                     \end{multlined}
  \end{align*}
\end{proposition}
\begin{proof}
  We start by estimating the approximation error. 
  By Lemma \ref{lemma:inputnorm},
  \begin{equation*}
    \|\matr(\Realiz(\Phiinputtwo (\axq))   )\|_2 \leq 1 -\delta.
  \end{equation*}
We temporarily restrict $\epsinv \in (0, \alpha/4)$.  Then, we have, for all $a\in \cD$,
  \begin{multline*}
    \|(\AtildeNbNq)^{-1} - \matr(\Realiz(\PhiinvA))(\axq) \|_2
    \\ 
    \begin{aligned}[t]
       \overset{\eqref{eq:PhiinvA-def}}&{=} \alpha\|(\alpha \AtildeNbNq)^{-1} - \matr(\Realiz(\Phiinv\sconc\Phiinputtwo))(\axq) \|_2
    \\ 
         \overset{\text{L.~\ref{lemma:Phiinputtwo}}}&{=} \alpha\|(\alpha \AtildeNbNq)^{-1} - \matr(\Realiz(\Phiinv)(-\alpha \vec(\AtildeNbNq) + \vec(\Id_{\nb })) \|_2
    \\ 
         \overset{\text{T.~\ref{th:inversion}}}&{\leq} \alpha \frac{\epsinv}{\alpha} = \epsinv.
    \end{aligned}
  \end{multline*}
  We now have to bound the depth and size of $\PhiinvA$. First, we remark that
  \begin{equation*}
m(\epsinv/\alpha,\delta) = \left\lceil \frac{\log\left(\Ccoer \epsinv /2\right)}{\log(1-\delta)}\right\rceil,
  \end{equation*}
  where $m(\cdot, \cdot)$ is defined in Theorem \ref{th:inversion}.
  Now, we use the fact that there exists $C_1 > 1$ such that, 
  for all $\epsinv \in (0,1)$,
  \begin{equation*}
  |\log(\Ccoer \epsinv/2)| \leq C_1(1+|\log\epsinv|).
  \end{equation*}
  Furthermore, there exists $C_2>1$ such that for all $\nb \in \N$, $\delta \geq C_2^{-1}$.
  Remark then that $\left|\log(1- y)\right| \geq y$ for all $y\in \left(0, 1\right)$, hence
    $\left|\log(1-\delta)\right|^{-1} \leq C_2 $.
  We infer that for all $\epsinv\in (0,1)$ 
  and
  for all $\nb \in\N$,
  \begin{equation*}
    m(\epsinv/\alpha, \delta) \leq C_1C_2 (2+|\log\epsinv|) \;.
  \end{equation*}
  Therefore, from Theorem \ref{th:inversion} we obtain that there exist
  constants $C_4, C_5>0$ dependent only on $\Ccoer$, $\Ccont$, and $d$
  such that
  \begin{equation*}
    \depth\left( \Phiinv \right) \leq C_4\left(1+\log(1+|\log\epsinv|) + \log(\nb ) \right) \cdot \left(1+ |\log\epsinv|  + \log(\nb ) +  \log(1+|\log\epsinv|) \right)
  \end{equation*}
  and 
  \begin{multline*}
    \size\left( \Phiinv \right)  \leq C_5 
(1+|\log\epsinv|)\nb ^{3} \big[  1+\log(1+|\log\epsinv|) + \log(\nb ) \big]^2 
\\ \times\big[ 1+|\log\epsinv| + \log(\nb ) +  \log(1+|\log\epsinv|)  \big].
  \end{multline*}
  We can now drop the restriction $\epsinv<\alpha/4$, adjusting the
    constants $C_4$ and $C_5$.
  Since, in addition, 
  \begin{equation*}
    \depth(\Phiinputtwo)=2,\qquad \size(\Phiinputtwo)\leq C_6\nb^3,
  \end{equation*}
  for $C_6 >0$ independent of $\nb $, we obtain the bounds on the depth and size
  of $\PhiinvA$.
\end{proof}
\subsubsection{Computation of the coefficients}
\begin{proposition}
  \label{prop:cu}
There exists a constant $\Ccu>0$ such that 
  for all $\nb \in \N$ and for all $\epsu\in (0,1)$, writing $\nq  = \nb+1$ 
  and 
  \begin{equation*}
    \Phicu \coloneqq \left( (\tbcf^\top \otimes \rev{(\AoneNbNq)^{-1/2}}, \bzero_{\nb }) \right) \sconc \PhiinvAepsu,
  \end{equation*}
  where $\CAinv$ is the constant from Lemma \ref{lemma:norminvA},
  we have
  \begin{equation*}
    \sup_{a\in\cD} \| \bcau - \Realiz(\Phicu)(\axq) \|_2 \leq \epsu
  \end{equation*}
  and
  \begin{align*}
    \depth(\Phicu)&
                    \begin{multlined}[t][.8\textwidth]
                      \leq \Ccu \big[1+\log(1+\left|\log\epsu\right| +\log\nb) + \log(\nb ) \big] \\ \times \big[1+ \left|\log\epsu\right|  + \log(\nb ) +  \log(1+\left|\log\epsu\right| + \log\nb) \big]
                    \end{multlined}
    \\
    \size(\Phicu)&
                     \begin{multlined}[t][.8\textwidth]
                     \leq\Ccu  \nb ^{3}
\big[1+\left|\log\epsu\right|\big] \big[ 1+ \log(1+\left|\log\epsu\right| +\log\nb) + \log(\nb )\big]^2 
\\ \times\big[1+ \left|\log\epsu\right| + \log(\nb ) +  \log(1+\left|\log\epsu\right|) \big]
   \;.
                     \end{multlined}
  \end{align*}
\end{proposition}
\begin{proof}
  For all $m, n, l\in \N$, let $\mathbf{B}\in \mathbb{R}^{m\times n}$ and
  $\mathbf{C}\in \mathbb{R}^{n\times l}$ \rev{and
    $\mathbf{X}\in\mathbb{R}^{m\times m}$}.
  Then
  \begin{equation*}
    \vec(\rev{\mathbf{X}}\mathbf{B}\mathbf{C}) = (\mathbf{C}^\top\otimes \rev{\mathbf{X}}) \vec(\mathbf{B}). 
  \end{equation*}
  This identity implies
  \begin{multline}
    \label{eq:cu-err-proof-1}
    \left(\tbcf^\top \otimes \rev{( \AoneNbNq )^{-1/2}} \right) \Realiz(\PhiinvAepsu)(\axq) 
  \\   = 
      \rev{( \AoneNbNq )^{-1/2}}\matr\left(\Realiz(\PhiinvAepsu)(\axq)\right) \tbcf.
  \end{multline}
We assume that $\nb^{3/2}\CAinv \|f\|_{L^2({Q})} \geq 1$; if this does not hold, 
  it
  is sufficient to temporarily restrict $\epsu < \nb^{3/2}\CAinv \|f\|_{L^2({Q})}$ and drop this
  restriction at the end of the proof by adjusting the constants. Therefore, for all $a\in \cD$,
  \begin{multline*}
    \| \bcau - \Realiz(\Phicu)(\axq)\|_2
    \\
    \begin{aligned}[t]
    \overset{\eqref{eq:cudef}}&{=}\| \rev{( \AoneNbNq )^{-1/2}}\left((\AtildeNbNq)^{-1} - \matr\left(\Realiz(\PhiinvAepsu)(\axq)\right) \right)\tbcf \|_2\\
    &\leq \rev{\|( \AoneNbNq )^{-1/2}\|_2}\|(\AtildeNbNq)^{-1} - \matr\left(\Realiz(\PhiinvAepsu)(\axq)\right) \|_2\|\tbcf \|_2\\
    \overset{\text{P.~\ref{prop:invA}}}& {\leq} \frac{\epsu}{\nb^{3/2}\CAinv\|f\|_{L^2({Q})}}\| (\AoneNbNq)^{-\rev{1/2}}\|^{\rev{2}}_2 \|\bcf\|_2\\
     \overset{\eqref{eq:bound-ctof}}&{\leq }
      \frac{\epsu}{\nb^{3/2}\CAinv\|f\|_{L^2({Q})}}\| (\AoneNbNq)^{-1}\|_2
            \sqrt{\nb}\|f\|_{L^2({Q})}\\
    \overset{\text{L.~\ref{lemma:norminvA}}}& {\leq }\frac{\epsu}{{\CAinv}\nb\|f\|_{L^2({Q})}}{\CAinv\nb}\|f\|_{L^2({Q})} = \epsu,
    \end{aligned}
  \end{multline*}
  where in the last three steps we have used 
  Proposition \ref{prop:invA}, bound \eqref{eq:bound-ctof}, and Lemma \ref{lemma:norminvA}.
   To derive the bounds on the size and
  depth of $\PhiinvAepsu$, first remark that 
  \begin{equation*}
    \|\tbcf^\top \otimes \rev{(\AoneNbNq)^{-1/2}} \|_0 \leq \nb ^{\rev{3}}.
  \end{equation*}
  Then, defining $\epsinv \coloneqq \epsu
    /(\nb^{3/2}\|f\|_{L^2({Q})}{\CAinv})$, there exists $C_1>0$ such that for all
    $\nb\in\N$ and for all $\epsu \in (0, 1)$,
    \begin{equation*}
      \left| \log\epsinv \right| \leq C_1\left( 1+ \left|  \log\epsu \right|  + \log\nb\right).
    \end{equation*}
    Inserting this bound in Proposition \ref{prop:invA}
    and applying Proposition \ref{prop:invA} concludes the proof.
\end{proof}
\subsection{Trunk network}
\label{sec:trunk}
The following emulation rates for the approximation of the polynomial basis are
a direct consequence of \cite[Proposition 2.13]{Opschoor2019}.
\begin{proposition}
  \label{prop:basis}
 There exists $\Cbasis>0$ such that, for all $\epsbasis\in (0, 1)$ and all $\nb \in \N$, 
 there exists a NN $\Phibasis$
 such that $\Realiz(\Phibasis) : \mathbb{R}^{d}\to \mathbb{R}^{\nb }$, 
 \begin{equation*}
   \max_{i \in \brange{\nb }} \| \phi_i - \left[\Realiz(\Phibasis)  \right]_i \|_{H^1({Q})} \leq \epsbasis,
 \end{equation*}
 and 
$$
\begin{array}{rcl}
   \depth(\Phibasis) &\leq& \Cbasis \left(1+ \left| \log\epsbasis \right|  + \nb ^{1/d}\right)\left( 1+ \log \nb\right) \\
   \size(\Phibasis)&\leq&
                     \Cbasis \left( \nb ^{2/d} + \nb ^{1/d}\left| \log\epsbasis \right| + \nb (1+\log \nb  + \left| \log\epsbasis \right|) \right).
\end{array}
$$
\end{proposition}
\subsection{Operator network expression rates}
\label{sec:ONetExpRt}
Combining the results from Sections \ref{sec:branch} and \ref{sec:trunk}, we obtain
the main result on the operator network
approximation of \eqref{eq:problem}. 
The structure of the operator network is schematically 
represented in Figures \ref{fig:ONet-scheme} and \ref{fig:ONet}.
\begin{theorem}
  \label{th:deepONet}
  There exists $C>0$ such that, for all $\epsilon \in (0,1)$, 
  for all $a\in\cD$ with $u^a=S(a)$ and solution operator $S$ 
  as defined in Section \ref{sec:problem},
  there exist
  \begin{enumerate}[label=(\alph*)]
  \item $\nb , \nq \in \N$,
  \item a set of points $\xenc\coloneqq \{\bx_1, \dots \bx_{\nq }\}\subset \overline{Q}$,
  \item two NNs $\Phibranch$ and $\Phitrunk$ with $\Realiz(\Phibranch) : \R^{\nq } \to \R^{\nb }$ and $\Realiz(\Phitrunk) : Q \to \R^{\nb }$,
  \end{enumerate}
such that
  \begin{enumerate}[label=(\roman*)]
  \item \label{item:NbNq}$\nb , \nq  \leq C (1+ \left| \log\epsilon \right|^d)$,
  \item \label{item:error} 
     the following error bound holds:
    \begin{equation*}
      \sup_{a\in\cD} \| u^a - (\cR_{\Phitrunk} \circ \cA_{\Phibranch}\circ\cE_{\xenc})(a)\|_{H^1({Q})}\leq \epsilon,
    \end{equation*}
  \item\label{item:size} 
    as $\epsilon \downarrow 0$,
    \begin{equation*}
      \depth(\Phibranch) = \cO\left(\left| \log\epsilon \right| ( \log\left| \log\epsilon \right| )\right),
      \qquad 
      \size(\Phibranch) = \cO\left(\left| \log\epsilon \right|^{3d+2} (\log\left| \log\epsilon \right|)^2\right),
    \end{equation*}
     and  
     \begin{equation*}
      \depth(\Phitrunk) = \cO\left(\left| \log\epsilon \right| ( \log\left| \log\epsilon \right| )\right),
      \qquad 
      \size(\Phitrunk) = \cO\left(\left| \log\epsilon \right|^{d+1}\right).
    \end{equation*}
  \end{enumerate}
\end{theorem}
\begin{proof}
Due to Lemma \ref{lemma:GNI}, 
there exist constants $\Cfe, \bfe, \Cq > 0$ such
that for all $\nb \in \N$, there exists $\nq  \leq \Cq\nb $ such that
\begin{equation*}
  \sup_{a\in \cD}\| u^a - \uaNbNq \|_{H^1({Q})} \leq \Cfe\exp(-\bfe \nb ^{1/d}),
\end{equation*}
where 
$\uaNbNq = \sum_{i=1}^{\nb } \left[\bcau  \right]_i \phi_i\in X_{\nb}$ 
is the
fully discrete Galerkin projection of $u^a$ onto $X_{\nb}$,
  \begin{equation*}
    \uaNbNq \in  X_{\nb}: \quad \baNq(\uaNbNq, v) = (f, v), \qquad \forall v\in X_{\nb}.
  \end{equation*}
  We assume, without loss of generality and for ease of notation, that
    $\Cfe\geq 1$.
Fix now 
\begin{equation}
  \label{eq:Nbfinal}
  p(\epsilon) =\left\lceil \frac{\left| \log(\epsilon/3) \right| + \log\Cfe}{\bfe} \right\rceil  + 1 , \qquad  \Nbeps= p(\epsilon)^d -1 \;, \qquad \Nqeps = \Nbeps + 1.
\end{equation}
Then,
we observe that \eqref{eq:Nbfinal} implies the existence of a constant $\Cb>0$ such that, for
all $\epsilon \in (0,1)$, we have $\Nbeps, \Nqeps \leq \Cb(1+\left| \log\epsilon \right|^d)$,
i.e., item \ref{item:NbNq} of the statement of the theorem.

We also define $\Cpol>0$ as a constant such that, for all $p\in \mathbb{N}$,
\begin{equation}
  \label{eq:poly-inverse}
  \|\nabla q\|_{L^2({Q})}\leq \Cpol p^2 \| q\|_{L^2({Q})}, \qquad \forall q\in\mathbb{Q}_p(Q).
\end{equation}
This inverse inequality follows straightforwardly from the classical Markov inequality in $(0,1)$,
with a tensorization argument (which yields that $\Cpol \sim \sqrt{d}$).
With $\Nbeps$ as in \eqref{eq:Nbfinal} define
  \begin{equation} \label{eq:alleps}
  \epsbasis \coloneqq \dfrac{\epsilon}{3\Nbeps(2 + \sup_{a\in \cD}\|u^a\|_{L^2({Q})})}, \quad 
  \epsu \coloneqq \dfrac{\epsilon}{3( 1 + \Cpol^2\Nbeps^{4/d} )^{1/2}\Nbeps^{1/2}},  \quad
  \epsfe \coloneqq \dfrac{\epsilon}{3}. 
\end{equation}
We assume that $\epsilon \in (0,1)$ implies $\epsbasis, \epsu \in (0,1)$. If
this is not the case, it is sufficient to restrict $\epsilon$ to values such
that $\epsbasis, \epsu\in (0,1)$; the restriction can be dropped at the end of
the proof. Due to \eqref{eq:Nbfinal},
\begin{equation}
  \label{eq:errorG}
  \sup_{a\in \cD}\| u^a - \uaNbNqeps \|_{H^1({Q})} \leq \epsfe.
\end{equation}
Define then
\begin{equation*}
  \Phibranch = \Phicufinal\qquad \text{and}\qquad \Phitrunk = \Phibasisfinal,
\end{equation*}
where the NNs $\Phicufinal$ and $\Phibasisfinal$ are defined in 
Propositions \ref{prop:cu} and \ref{prop:basis}, respectively.
\paragraph{\textbf{Error estimate}}
For all $a\in \cD$, 
\begin{multline*}
  \| u^a - \left( \Realiz(\Phibranch) (\axq) \right)\cdot \Realiz(\Phitrunk) \|_{H^1({Q})}\\
  \leq 
  \| u^a - \uaNbNqeps \|_{H^1({Q})} + 
  \| \uaNbNqeps - \left( \Realiz(\Phibranch) (\axq) \right)\cdot \Realiz(\Phitrunk)
  \|_{H^1({Q})}
  \eqqcolon (I) + (II).
\end{multline*}
We have already established that $(I)\leq \epsfe  = \epsilon/3$. 

Consider term $(II)$.  We have
\begin{align*}
  (II) & =   \| \sum_{i=1}^{\Nbeps} \left(  \left[ \bcaueps \right]_i \phi_i  -
         \left[ \Realiz(\Phicufinal) (\axq)\right]_i \left[ \Realiz(\Phibasisfinal) \right]_i \right)\|_{H^1({Q})}
  \\ & \leq
       \begin{multlined}[t][.8\textwidth]
       \|\sum_{i=1}^{\Nbeps} \left(  \left[ \bcaueps \right]_i -
         \left[ \Realiz(\Phicufinal) (\axq)\right]_i \right) \phi_i\|_{H^1({Q})}
       \\ + \|\sum_{i=1}^{\Nbeps}\left[ \Realiz(\Phicufinal) \right]_i \left(  \phi_i  -
          \left[ \Realiz(\Phibasisfinal) \right]_i \right)\|_{H^1({Q})}
       \end{multlined}
       \\ & \eqqcolon (IIa) + (IIb).
\end{align*}
Denote, for all $i\in\brange{\Nbeps}$,
\begin{equation*}
  \eta_i \coloneqq \left[ \bcaueps \right]_i - \left[ \Realiz(\Phicufinal) (\axq)\right]_i.
\end{equation*}
Using 
the Cauchy-Schwarz inequality, the bound $\|\phi_i\|_{L^2(Q)}\leq 1$,
the polynomial inverse inequality \eqref{eq:poly-inverse}
and Proposition \ref{prop:cu}, 
we obtain
\begin{align*}
  (IIa)^2&\leq \| \sum_{i=1}^{\Nbeps}\eta_i \phi_i\|^2_{H^1({{Q}})}
           = \int_{Q}  \left( \sum_{i=1}^{\Nbeps}\eta_i\phi_i \right)^2 
           + \int_{Q}\left( \sum_{i=1}^{\Nbeps}\eta_i\nabla\phi_i \right)^2 
     \\ 
  \overset{\text{C-S}}&{\leq} \left(\sum_{i=1}^{\Nbeps}  \eta_i^2 \right)
  \sum_{i=1}^{\Nbeps}\left( \|\phi_i\|^2_{L^2(Q)} +\|\nabla\phi_i\|_{L^2(Q)}^2 \right)
  \\ 
       \overset{\eqref{eq:poly-inverse}}&{\leq} 
        \| \bcaueps - \Realiz(\Phicufinal) (\axq)\|_2^2\left( 1 + \Cpol^2\Nbeps^{4/d} \right)\Nbeps
  \\ 
       \overset{\text{P.~\ref{prop:cu}}}&{\leq} \epsu^2 \left( 1 + \Cpol^2\Nbeps^{4/d} \right)\Nbeps
  \\
         \overset{\eqref{eq:alleps}}&{\leq} \left( \frac{\epsilon}{3} \right)^2.
\end{align*}
Next, we estimate
\begin{equation}
  \label{eq:Phicunorm}
\begin{aligned}
  \|\Realiz(\Phicufinal) (\axq) \|_2
  &\leq  
    \|\Realiz(\Phicufinal) (\axq) - \bcaueps \|_2 + \|\bcaueps\|_2 
 \\   \overset{\text{P.~\ref{prop:cu}},\eqref{eq:L2-l2}}&{\leq} 1 + \CLell\Nbeps^{1/2}\| \uaNbNqeps \|_{L^2({Q})} \\
  &\leq 1 + \CLell\Nbeps^{1/2}\| \uaNbNqeps - u^a\|_{L^2({Q})} + \CLell\Nbeps^{1/2}\| u^a \|_{L^2({Q})}\\
  \overset{\eqref{eq:errorG}}&{\leq} \left(2 + \| u^a \|_{L^2({Q})}  \right)\CLell\Nbeps^{1/2}.
\end{aligned}
\end{equation}
Then,
\begin{align*}
  (IIb)^2 \overset{\text{C-S}}&{\leq}
  \|\Realiz(\Phicufinal) (\axq) \|^2_2 
            \sum_{i=1}^{\Nbeps}\|\phi_i  - \left[ \Realiz(\Phibasisfinal) \right]_i \|^2_{H^1({Q})}
  \\
  &\leq\Nbeps \|\Realiz(\Phicufinal) (\axq) \|^2_2 \max_{i\in\brange{\Nbeps}}
            \|\phi_i  - \left[ \Realiz(\Phibasisfinal) \right]_i \|^2_{H^1({Q})}\\
  \overset{\text{P.~\ref{prop:basis}}, \eqref{eq:Phicunorm}}&{\leq} \Nbeps(2+\|u_a\|_{L^2({Q})})^2 \CLell^2\Nbeps\epsbasis^2
       \\ \overset{\eqref{eq:alleps}}&{\leq} \left( \frac{\epsilon}{3}  \right)^2.
\end{align*}
We can conclude that
\begin{equation*}
  \| u^a - \left( \Realiz(\Phibranch) (\axq) \right)\cdot \Realiz(\Phitrunk) \|_{H^1({Q})} 
  \leq (I) + (IIa) + (IIb) \leq \epsilon.
\end{equation*}
\paragraph{\textbf{Depth and size bounds}}
Using \eqref{eq:Nbfinal} and the definitions \eqref{eq:alleps}, we obtain that there
exists a constant $C_1 >0$ such that,
for all $\epsilon \in (0, 1)$, 
\begin{equation*}
  1 + \max\left( \left|\log\epsbasis\right|, \left|\log\epsfe\right|, \left|\log\epsu\right| \right)\leq C_1 (1+ \left|\log\epsilon\right|).
\end{equation*}
We infer then, 
from Proposition \ref{prop:cu}, that there exists $C_2 >0$ such that, 
for all $\epsilon \in (0, 1)$,
\begin{equation*}
  \depth(\Phicufinal) \leq C_2\left( 1+\log(1+\left| \log\epsilon \right|^d) \right) \left( 1+\left| \log\epsilon \right| \right)
\end{equation*}
and 
\begin{equation*}
  \size(\Phicufinal) \leq C_2\left( 1+\left| \log\epsilon \right|^d \right)^{3} \left( 1+\left| \log\epsilon \right| \right)^2 \left( 1+ \log(1+\left| \log\epsilon \right|^d) \right)^2.
\end{equation*}
Furthermore, from Proposition \ref{prop:basis}, we have that there exists
$C_3>0$ such that for all $\epsilon\in (0,1)$
\begin{equation*}
  \depth(\Phibasisfinal) \leq C_3 (1+\left| \log\epsilon \right|)\left(1+\log(1+\left| \log\epsilon \right|^d)  \right)
\end{equation*}
and 
\begin{equation*}
  \size(\Phibasisfinal) \leq C_3 (1+\left| \log\epsilon \right|^{d+1}).
\end{equation*}
Using the definition of $\Phitrunk$ and $\Phibranch$ gives Item \ref{item:size}
and concludes the proof.
  \end{proof}
  \begin{remark}
  The implicit constants in the size bounds of the operator 
  networks in Theorem \ref{th:deepONet} and in the theorems of the upcoming sections depend, 
  in general, exponentially on the dimension $d$.
  \end{remark}

  \begin{remark}
      \label{remark:fmap}
      In Theorem \ref{th:deepONet}, we have considered the data-to-solution 
      map $a\mapsto u^a$ for a fixed, given right-hand side $f$. 
      The present analysis may be extended to the operator
$S: (a,c, f) \mapsto u$
      where $u$ is the solution to the reaction-diffusion equation
$$
  -\nabla \cdot (a\nabla u)  + c u = f \; \text{ in } \; \Omega
$$
with 
      (analytic in $\overline{Q}$ and $Q$-periodic) 
      positive diffusion coefficient function $a$,
      non negative reaction coefficient function $c$ and source term $f$. 
      This requires a straightforward
modification of the branch network so that
      \begin{enumerate}
      \item it takes the point evaluations of $f$ at $\xenc$ as input and outputs 
        an exponentially consistent numerical quadrature approximation of $\bcf$ in \eqref{eq:bcf-def};
      \item the approximation of $\bcf$ 
        is then passed to a network approximating
        matrix-vector multiplication (see \cite[Proposition 3.7]{Kutyniok2019})
        with the output of the network of Proposition \ref{prop:invA}.
      \end{enumerate}
      The construction of the remaining parts of the operator network 
      follows along the same lines.
    \end{remark}
  \subsection{Parametric diffusion coefficient}
  \label{sec:parametric}
In many applications, for example in uncertainty quantification, one is interested in
    the case where the diffusion coefficient in \eqref{eq:problem} is parametric. 
This is naturally accommodated for by composition with 
     solution operator networks and we briefly detail this here.
  Specifically, suppose that there exists
  $\dimp\in \N$ and a compact parameter set $\cP\subset\R^{\dimp}$ 
  such that $\ap : \cP\to \Hol(\Omega)$
  and that there exist constants $\amin, \Cp, \bp, \alphap, \Ap, A_\psi > 0$,  
  and functions $\psi_i : Q\to \R$ 
  and $a_i:\cP \to \R$,  $i\in \N$, 
  such that
  \begin{equation}
    \label{eq:amin-par}
    \inf_{\by\in\cP}\inf_{\bx\in{Q}}\ap(\by)(\bx)\geq \amin,
  \end{equation}
  that
  \begin{equation}
    \label{eq:a-decomp}
    \forall \np\in \N,\qquad
    \sup_{\by\in \cP} \|\ap(\by) - \sum_{i=1}^{\np}a_i(\by) \psi_i\|_{L^\infty({Q})} 
     \leq \Cp\exp(-\bp\np^{\alphap}),
  \end{equation}
with
  \begin{equation} \label{eq:apsi-bound}
    \forall i\in \N,\qquad 
   \psi_i\in \Hol(\Omega; A_\psi), \qquad a_i\in \Hol(\cP; \Ap),
  \end{equation}
  and that
  \begin{equation}
    \label{eq:apsi-bound-2}
    \sup_{\by\in\cP} \sum_{i=1}^{\infty} |a_i(\by)|\leq \Ap.
  \end{equation}
  Here, we use the same constant $\Ap$ in the second hypothesis in
    \eqref{eq:apsi-bound} and in \eqref{eq:apsi-bound-2} only to simplify notation.
  For all $\by \in \cP$, we denote $u_{\by}\in X$ the solution to
  \begin{equation}
    \label{eq:parametric}
    -\nabla\cdot(\ap(\by)\nabla u_{\by} ) = f, \qquad \text{in }\Omega.
  \end{equation}
\begin{remark}
    Diffusion coefficient functions that can be written in Fourier series as
    \begin{equation*}
      \ap(\by)(\bx) = \sum_{\bk\in \Z^d }(a_{\bk}(\by) +i b_{\bk}(\by))e^{i2\pi \bk \cdot \bx},
    \end{equation*}
    where $a_{\bk}, b_{\bk}\in \Hol(\cP, \Ap)$ are chosen so
      that $\ap$ is a real function for all $\by$, with exponential decrease of
    $\sup_{\by\in\cP}(|a_{\bk}(\by)| + |b_{\bk}(\by)|)$ with respect to $|\bk|$, and such that $\ap$ is
    uniformly bounded from below by a positive constant 
    in the sense that \eqref{eq:amin-par} holds for its real part,
    fulfill conditions \eqref{eq:a-decomp}, \eqref{eq:apsi-bound} and \eqref{eq:apsi-bound-2}.
  \end{remark}

  \begin{lemma}
    \label{lemma:ai-approx}
    There exists $C>0$ such that for all $\np\in\N$ and for all $\epsilon\in(0,1)$, 
    there exists a NN $\Phiaone$ with 
    input dimension $\dimp$ and
    output dimension $\np$ such that
    \begin{equation}
    \label{eq:ai-approx}
      \max_{i=1, \dots, \np}\|a_i - \left[ \Realiz(\Phiaone) \right]_i\|_{L^\infty(\cP)} \leq \epsilon
    \end{equation}
    and that $\depth(\Phiaone)\leq C(1+\left| \log\epsilon \right|)(1+\log\left| \log\epsilon \right|)$ and
    $\size(\Phiaone)\leq C(1+\left| \log\epsilon \right|^{\dimp+1})\np$.
  \end{lemma}
  \begin{proof}
    The statement follows from a parallelization of the network of \cite[Theorem 3.6]{Opschoor2019}.
      \end{proof}
  \begin{theorem}
\label{th:parametric} 
   Let $\dimp\in \N$ and let $\ap$ and $u_{\by}$ be defined as above. There exists $C>0$ such that, for all $\epsilon \in (0,1)$, there exist
  \begin{enumerate}[label=(\alph*)]
  \item $\nb \in \N$,
  \item two NNs $\Phibranch$ and $\Phitrunk$ with $\Realiz(\Phibranch) : \R^{\dimp } \to \R^{\nb }$ and $\Realiz(\Phitrunk) : Q \to \R^{\nb }$,
  \end{enumerate}
such that
  \begin{enumerate}[label=(\roman*)]
  \item \label{item:par-NbNq}$\nb \leq C (1+ \left| \log\epsilon \right|^d)$,
  \item \label{item:par-error}  
        the following error estimate holds:
    \begin{equation*}
      \sup_{\by\in\cP} \| u_{\by} - \left( \Realiz(\Phibranch) (\by) \right)\cdot \Realiz(\Phitrunk) \|_{H^1({Q})}\leq \epsilon,
    \end{equation*}
  \item\label{item:par-size} as $\epsilon \downarrow 0$,
    \begin{align*}
      \depth(\Phibranch) &= \cO\left(\left| \log\epsilon \right| ( \log\left| \log\epsilon \right| )\right),
      \\
      \size(\Phibranch) &= \cO\left(\left| \log\epsilon \right|^{3d+2} (\log\left| \log\epsilon \right|)^2 + \left| \log\epsilon \right|^{1+\dimp + 1/\alphap}\right),
    \end{align*}
     and  
     \begin{equation*}
      \depth(\Phitrunk) = \cO\left(\left| \log\epsilon \right| ( \log\left| \log\epsilon \right| )\right),
      \qquad 
      \size(\Phitrunk) = \cO\left(\left| \log\epsilon \right|^{d+1}\right).
    \end{equation*}
  \end{enumerate}
\end{theorem}
\begin{proof}
The proof proceeds in several steps. 
    We first prove a consistency bound,
    then detail the construction of the ONet, 
    and conclude with verification of the asserted bounds on the depth and size of the ONet.

  Let $C_L>0$ be the constant such that, given $\ap_1, \ap_2 \in L^\infty(\Omega)$
  such that
  \begin{equation}
    \label{eq:abounds-extended}
    0 < \frac{\amin}{4}\leq \ap_i \leq \max(\amax, (1+\Ap)A_\psi), \qquad 
           \text{a.~e.~in }{Q} \text{ and for }i=1,2,
  \end{equation}
  and $u_i =S(\ap_i)$, $i=1,2$,
  then
  \begin{equation*}
   \| u_1 - u_2\|_{H^1({Q})} \leq C_L \|\ap_1-\ap_2\|_{L^\infty({Q})},
  \end{equation*}
  see Lemma \ref{lemma:lip-sol}.
  We suppose, without loss of generality and for ease of notation, that
  $A_\psi\geq 1$ and $C_L\geq 1$.
  Let now $\np $ be the smallest integer such that
  \begin{equation}
    \label{eq:npeps}
    C_L\Cp\exp(-\bp \np^{\alphap}) \leq \min\left(\frac{\epsilon}{3}, \frac{\amin}{2}\right).
  \end{equation}
  This implies that there exists a constant $C_1 > 0$ 
  (depending only on $C_L, \Cp, \bp, \amin$) 
  such that 
  \begin{equation*}
    \np \leq C_1(1+\left| \log\epsilon \right|^{1/\alphap})
  \end{equation*}
  and that, due to 
    \eqref{eq:amin-par}, \eqref{eq:a-decomp},
    \eqref{eq:apsi-bound}, and \eqref{eq:apsi-bound-2},
  \begin{equation}
    \label{eq:asum-bound}
    \inf_{\by\in\cP}\inf_{x\in{Q}}\sum_{i=1}^{\np}a_i(\by)\psi_i (x)\geq \frac{\amin}{2}, \qquad 
    \sup_{\by\in\cP}\|\sum_{i=1}^{\np}a_i(\by)\psi_i \|_{L^\infty({Q})}\leq \Ap A_\psi.
  \end{equation}
  Let also
  \begin{equation}
    \label{eq:epsp-tap}
    \epsp \coloneqq \frac{1}{\np A_\psi}\min\left(\frac{\epsilon}{3C_L}, \frac{\amin}{4}  \right), 
     \qquad \tap \coloneqq \sum_{i=1}^{\np} \left[\Realiz(\Phiaoneepsp)  \right]_i\psi_i,
  \end{equation}
  where the network $\Phiaoneepsp$ is defined in Lemma \ref{lemma:ai-approx}.
  We now show that $\tap$ fulfills conditions like \eqref{eq:abounds} and
  \eqref{eq:abound} (with updated values of the constants $\amin$, $\amax$,
  $A_{\cD}$), uniformly with respect to $\np$ and $\epsp$.
  From \eqref{eq:ai-approx}
  and \eqref{eq:asum-bound}, 
  it follows that, for all $k\in\N_0$,
  \begin{equation}
    \label{eq:tap-analytic}
    \begin{aligned}
      \sup_{\by\in\cP}\|\tap (\by)\|_{W^{k, \infty}({Q})}
      & \leq 
      \sup_{\by\in \cP}  \sum_{i=1}^{\np} \left(| a_i(\by) - \left[ \Realiz(\Phiaoneepsp)(\by) \right]_i| + |a_i(\by)|\right) \| \psi_i\|_{W^{k, \infty}({Q})}
      \\ & \leq 
      \left( \np \epsp   + \Ap\right) A_{\psi}^{k+1}k!
      \leq \left( 1  + \Ap\right) A_{\psi}^{k+1}k!.
    \end{aligned}
  \end{equation}
  Furthermore, for all $\by\in\cP$ and all $x\in{Q}$,
  \begin{equation}
    \label{eq:tap-lower}
    \begin{aligned}
      \tap(\by)(x)
      & \geq 
     \sum_{i=1}^{\np} \left( \left[ \Realiz(\Phiaoneepsp) (\by)\right]_i - a_i(\by) \right) \psi_i(x)+ \sum_{i=1}^{\np}a_i(\by)\psi_i(x)
     \\ & \geq
     \frac{\amin}{2} - \np\epsp A_\psi,
     \\ & \geq \frac{\amin}{4}.
    \end{aligned}
  \end{equation}
  Here we have used \eqref{eq:apsi-bound-2}, \eqref{eq:ai-approx}, \eqref{eq:asum-bound},
  and the definition of $\epsp $ in \eqref{eq:epsp-tap}.

  \paragraph{\textbf{Construction of the operator network and error estimate}}
  For $q\in \N$, $\nq = q^d$, 
  we introduce the matrix $\Vpsi \in \R^{\nq\times\np}$ 
  with entries 
  \begin{equation}
    \label{eq:Vpsi-def}
    \left[ \Vpsi \right]_{ij} = \psi_{j}(\xq_{i}), \qquad i=1,\dots, \nq ,\; j=1, \dots, \np,
  \end{equation}
  where $\xq_1, \dots, \xq_{\nq}$ are the quadrature nodes introduced in Section \ref{sec:pol-quad}.
  Then the NN
  \begin{equation}
    \label{eq:Phitap}
    \Phitap = \left( \left( \Vpsi, \bzero_{\nq} \right) \right) \sconc \Phiaoneepsp
  \end{equation}
  has realization such that
  \begin{equation*}
    \Realiz(\Phitap)(\by) =
    \begin{pmatrix}
      \tap(\by)(\xq_1)\\
      \vdots\\
      \tap(\by)(\xq_{\nq})
    \end{pmatrix}.
  \end{equation*}
  Let $\tu_{\by}\in X$ denote, for each $\by\in \cP$, 
  the solution to 
  \begin{equation*}
    - \nabla\cdot(\tap(\by)\nabla\tu_{\by}) = f \quad \text{in }\Omega.
  \end{equation*}
  Thanks to \eqref{eq:tap-analytic}, \eqref{eq:tap-lower}, and to Theorem
  \ref{th:deepONet}, there exists a constant $C_2$ independent of $\epsilon$,
  $\nq\in \N$ such that $\nq\leq C_2(1+\left| \log\epsilon \right|)$, 
  and networks $\tPhibranch$ and $\tPhitrunk$ such that
  \begin{equation*}
   \forall \by \in \cP: \quad 
    \| \tu_{\by} -\left(\Realiz(\tPhibranch) \circ\Realiz(\Phitap)  \right) (\by)\cdot\Realiz(\tPhitrunk) \|_{H^1({Q})} 
    \leq \frac{\epsilon}{3}.
  \end{equation*}
  Furthermore, for all
  $\by\in \cP$, $\ap(\by)$ and $\tap(\by)$ satisfy the conditions in
  \eqref{eq:abounds-extended}, hence for all $\by \in \cP$
  \begin{align*}
    \| u_{\by} - \tu_{\by} \|_{H^1({Q})}
    & \leq
      C_L \| \ap(\by) - \tap(\by) \|_{L^\infty({Q})}
    \\ & \leq 
          C_L \left(
          \| \ap(\by) - \sum_{i=1}^{\np}a_i(\by)\psi_i \|_{L^\infty({Q})}
                      +\| \sum_{i=1}^{\np}\left( a_i(\by) - \left[ \Realiz(\Phiaoneepsp) (\by) \right]_i\right)\psi_i \|_{L^\infty({Q})}
 \right)
    \\  &\leq 
         \frac{\epsilon}{3}
          +C_L\np \epsp A_\psi
     \\ &\leq  \frac{2}{3}\epsilon,
  \end{align*}
  where we have used \eqref{eq:a-decomp}, \eqref{eq:ai-approx}, \eqref{eq:npeps}, and
    \eqref{eq:epsp-tap} in the third inequality and \eqref{eq:epsp-tap} in the
    last one.
  We deduce that
  \begin{equation*}
    \sup_{\by\in \cP}\| u_{\by} - \left(\Realiz(\tPhibranch) \circ\Realiz(\Phitap)  \right) (\by)\cdot\Realiz(\tPhitrunk) \|_{H^1({Q})} \leq \epsilon,
  \end{equation*}
  which is Item \ref{item:par-error}, with
  \begin{equation*}
    \Phibranch \coloneqq \tPhibranch \sconc \Phitap, \qquad
    \Phitrunk \coloneqq \tPhitrunk.
  \end{equation*}
  \paragraph{\textbf{Depth and size bounds}}
  The bounds on the depth and size of $\Phitrunk$ can be inferred directly from
  Theorem \ref{th:deepONet}.
  To compute bounds on the size and depth of $\Phibranch$, note that, by Lemma \ref{lemma:ai-approx},
  there exist $C_3, C_4, C_5, C_6$ independent of $\epsilon$ such that
  \begin{equation}
    \label{eq:depth-Phia}
  \begin{aligned}
    \depth(\Phiaoneepsp)
   & \leq C_3(1+\left| \log\epsp  \right|)(1+\log\left| \log\epsp  \right|)
   \leq
   C_4(1+\left| \log\epsilon \right|)(1+\log\left| \log\epsilon \right|)
  \end{aligned}
  \end{equation}
  and 
  \begin{equation}
    \label{eq:size-Phia}
  \begin{aligned}
    \size(\Phiaoneepsp)
   & \leq C_5(1+\left| \log\epsp  \right|^{\dimp+1})\np
   \leq
   C_6(1+\left| \log\epsilon \right|^{\dimp+1+1/\alphap}).
  \end{aligned}
  \end{equation}
  Furthermore, there exists $C_7$ independent of $\epsilon$ such that
  \begin{equation}
    \label{eq:sizeV}
    \| \Vpsi \|_0 \leq \np\nq \leq C_7(1+\left| \log\epsilon \right|^{1+1/\alphap}).
  \end{equation}
  From \eqref{eq:size-Phia} and \eqref{eq:sizeV} it follows that
  \begin{equation}
   \label{eq:depthsizePhitap} 
   \depth(\Phitap) \leq C_8 (1+\left| \log\epsilon \right|)(1+\log\left| \log\epsilon \right|), \qquad
   \size(\Phitap) \leq C_9 (1+\left| \log\epsilon \right|^{\dimp+1+1/\alphap}),
  \end{equation}
  for constants $C_8, C_9$ independent of $\epsilon$.
  Combining the bounds in \eqref{eq:depthsizePhitap} with the bounds on the
  depth and size of $\tPhibranch$ coming from Theorem \ref{th:deepONet}
  concludes the proof.
\end{proof}
  \begin{remark}
    If each function $a_i$ does not depend on all the parameters but only on a
    subset of them, the size bound of Theorem \ref{th:parametric} results in an
    overestimation. Specifically, for all $i\in\N$, let $\cP_i$ be the domain of
    $a_i$ and denote $\dpi \coloneqq \dim(\cP_i)$.
    Then, 
    modifying Lemma \ref{lemma:ai-approx} so that the subnetworks
      approximating each of the $a_i$ only take a $\dpi$-dimensional input,
    we obtain
    in Theorem \ref{th:parametric} that there exists a constant
    $c>0$ independent of $\epsilon$ such that for
    $\epsilon \downarrow 0$,
    \begin{equation*}
      \size(\Phibranch)
      = \cO\left(\left| \log\epsilon \right|^{3d+2} (\log\left| \log\epsilon \right|)^2
        + \sum_{i=1}^{c\left| \log\epsilon \right|^{1/\alphap}}\left| \log\epsilon \right|^{1+\dpi }\right).
    \end{equation*}
    Clearly, setting $\dpi = \dimp$ for all $i$ in the equation above gives the
    estimate in Theorem \ref{th:parametric}.
  \end{remark}
  \begin{remark}
    \label{remark:kutyniok}
    Similar results to Theorem \ref{th:parametric} can be obtained through the
    technique in \cite{Kutyniok2019}, by using the exponential convergence of
    polynomial approximations to the functions in the solution manifold
    $\mathcal{M} = \{u(\by): \by\in\cP\}$ to derive an upper bound on the
    $n$-width of $\mathcal{M}$.
  \end{remark}
\section{Generalizations}
\label{sec:generalizations}
All steps of the analysis of ONet emulation rates for the 
coefficient-to-solution map of \eqref{eq:problem} 
directly generalize to other, structurally similar, linear 
divergence-form elliptic PDEs. 
We illustrate the extension of the preceding result by two of these: 
anisotropic diffusion-reaction equations and linear elastostatics.
\subsection{Linear anisotropic diffusion-reaction equations}
\label{sec:ADR}
\subsubsection{Definition of the problem}
\label{sec:ADR-prob}
We consider again the torus $\Omega = ( \R/\Z )^d$.
For a constant $A_{\tcD}>0$, 
    introduce the set of admissible data 
\begin{equation*}
  \tcD \subset \Hol(\Omega; A_{\tcD})^{d\times d}\times \Hol(\Omega; A_{\tcD})
\end{equation*}
of pairs $(\argsADR)$ and suppose 
there exist $Q_0\subset Q$, $\amin, \cmin>0$ such that 
for all $(\bA, c)\in\tcD$,
\begin{itemize}
\item $\bA$ is symmetric and is uniformly positive definite, i.e., $\bA_{ij} =
  \bA_{ji}$ and
\begin{equation*}
\forall \bx\in {Q},\; \forall \bxi\in \mathbb{R}^d,\quad 
\bxi^\top \bA(\bx) \bxi \geq \amin |\bxi|^2 \;,
\end{equation*}
\item $c(\bx)\geq $ for all $\bx\in Q_0$.
\end{itemize}
For all $(\argsADR )\in \tcD$, the bilinear form
$\frab^{(\argsADR)}(\cdot,\cdot): H^1(\Omega) \times H^1(\Omega) =\Hper^1(Q)
\times \Hper^1(Q)\to \mathbb{R}$
given by 
\begin{equation*}
\bADR(w,v)
\coloneqq
\int_{Q} \left( (\bA\nabla w)\cdot \nabla v +  c w v \right)
\end{equation*}
is coercive, i.e., 
there exists a constant $\alpha_0 > 0$ independent of $(\argsADR)$ 
such that 
\begin{equation*}
\forall v\in H^1(\Omega),\quad 
\bADR(v,v) \geq  \alpha_0 \| v \|_{H^1({Q})}^2 \;.
\end{equation*}
The continuity of the form 
$\bADR(\cdot, \cdot)$ on 
$H^1(\Omega) \times H^1(\Omega)\to \mathbb{R}$
being evident,
the Lax-Milgram Lemma implies that
for every $f\in \Hol(\Omega)$ there exists a unique solution 
$$
u\in H^1(\Omega):\quad 
\bADR(u,v) = (f,v) \quad \forall v\in H^1(\Omega)=\Hper^1(Q).
$$
For given, fixed $f\in \Hol(\Omega)$, the 
coefficient-to-solution map
$$
\SADR :  (\argsADR) \mapsto u
$$
is analytic. 
Furthermore, there exists $A_{\tcU}>0$ such that
  \begin{equation*}
  \SADR(\tcD)\subset \Hol(\Omega;
A_{\tcU}),
  \end{equation*}
which can be proven as in Lemma \ref{lemma:regularity}.
\subsubsection{Operator network approximation}
\label{sec:OpNetAppr}
We introduce, for all $\nq\in\N$ such that $q \coloneqq \nq^{1/d}\in \N$, the encoding
operator
$\cEADR: C(\Omega)^{d\times d} \times C(\Omega)\to \R^{d^2\nq +\nq}$ such that
\begin{equation*}
  \cEADR(\argsADR) =
  \begin{pmatrix}
    \vec(\bA(\xq_1))\\
    \vdots\\
    \vec(\bA(\xq_{\nq}))\\
    c(\xq_1)\\
    \vdots\\
    c(\xq_{\nq})
  \end{pmatrix},
\end{equation*}
where $\xenc = \xq_1, \dots, \xq_{\nq}$ are the points from Section \ref{sec:pol-quad}.
Theorem \ref{th:deepONet} can then be extended to this class of
reaction-diffusion equations.
\begin{theorem} \label{th:adr}
Theorem \ref{th:deepONet} holds with $a\in\cD$ replaced by 
$(\argsADR)\in\tcD$, $S(a)$ replaced by $\SADR(\argsADR)$, and $\cE_{\xenc}(a)$
replaced by $\cEADR(\argsADR)$.
\end{theorem}

For the proof, for all $p\in \N$, writing $\tnb = p^d$, 
consider the basis functions 
$\{\tphi_i\}_{i=1}^{\tnb }$ of $\Q_p(Q)$ such that
\begin{equation*}
  \{\tphi_i \}_{i=1}^{\tnb} = \{\phi_i\}_{i=1}^{\nb} \cup \{L_0 \otimes \dots \otimes L_0\},
\end{equation*}
where $L_0\equiv 1$ is the constant, unit function in $(0,1)$.
Define further
  \begin{equation*}
    \tX_{\nb}\coloneqq \spn(\{\tphi_1, \dots, \tphi_{\tnb}\})\subset \tX \coloneqq H^1(\Omega)=\Hper^1(Q).
  \end{equation*}

In order to prove Theorem \ref{th:adr}, we have to replace the input layer
network introduced in Lemma \ref{lemma:Phiinputone} with an input layer adapted
for anisotropic diffusion-reaction problems, as introduced in Lemma
\ref{lemma:adr-input} below.
For $k\in \brange{\nq}$, we introduce $\tbD(\xkq)$ such that
\begin{equation*}
 \tbD^{ij}_{mn} (\xkq)= \wkq (\partial_{x_n}\tphi_j) (\xkq)(\partial_{x_m}\tphi_i)(\xkq), \qquad (i,j)\in\brange{\tnb}^2, \, (m, n)\in\brange{d}^2.
\end{equation*}
Furthermore, let $\frav : \N^2 \to \N$ be the reordering such that for any matrix
$\bA$, 
\begin{equation}
  \label{eq:frav-1}
\vec(\bA)_{\frav(i, j)} = \bA_{ij}. 
\end{equation}
We introduce the operation
$\tvec:\R^{\tnb\times\tnb\times d\times d}\to \R^{\tnb^2\times d^2}$
\begin{equation*}
  \tvec(\tbD(\xkq))_{\frav(i,j)\frav(m,n)} = \tbD^{ij}_{mn}(\xkq).
\end{equation*}
Finally, define 
\begin{equation*}
  \hbM_{ij}(\xkq)   = \wkq \tphi_i(\xkq) \tphi_j(\xkq).
\end{equation*}
\begin{lemma}
  \label{lemma:adr-input}
For all $\alpha \in \mathbb{R}$ the one-layer NN
\begin{equation*}
  \Phiinputadr
  \coloneqq \left(
    \left(-\alpha \left[\tvec(\tbD(\xq_1)) | \dots |\tvec(\tbD(\xq_{\nq }))| \vec(\hbM(\xq_1) | \dots
      | \vec(\hbM(\xq_{\nq}))\right], \bzero_{\tnb^2}\right)
  \right)
\end{equation*}
is such that 
\begin{equation*}
  \matr\left(\Realiz(\Phiinputadr) (\cEADR(\argsADR))   \right)_{ij} = 
  -\alpha \frab^{(\argsADR)} (\tphi_j, \tphi_i)
\end{equation*}
and $\size(\Phiinputadr) \leq (d^2+1)\tnb ^2\nq  $ .
\end{lemma}
\begin{proof}
  We have
  \begin{multline*}
    \left[\Realiz\left( \Phiinputadr \right) (\cEADR(\argsADR))   \right]_{\frav(i,j)}
    \\
    \begin{aligned}[t]
    & = -\alpha\sum_{k=1}^{\nq }\wkq \left( \sum_{m, n=1}^d\left[ \bA_{mn}(\partial_{x_n}\tphi_j)(\partial_{x_m}\tphi_i)\right](\xkq)  + \left[ c \tphi_j\tphi_i \right] (\xkq)\right)\\
    &= -\alpha\sum_{k=1}^{\nq }\wkq \left( \left[ (\bA\nabla \tphi_i)\cdot (\nabla\tphi_j) \right](\xkq)  + \left[ c \tphi_j\tphi_i \right] (\xkq)\right)
    ,
    \end{aligned}
  \end{multline*}
  hence the equality after matricization. The size bound follows from the fact
  that
  \begin{equation*}
    \| \tbD(\xkq) \|_0 \leq d^2\tnb ^2,\qquad \| \hbM (\xkq) \|_0 \leq \tnb^2,
  \end{equation*}
  for all $k\in \brange{\nq }$.
\end{proof}

We can now prove Theorem \ref{th:adr}.

\begin{proof}[Proof of Theorem \ref{th:adr}]
  The proof follows along the same lines as the proof of Theorem \ref{th:deepONet}.
  In particular, in the construction of $\Phibranch$, the input network $\Phiinputone$ and Lemma
  \ref{lemma:Phiinputone} are replaced by the network
  $\Phiinputadr$ and Lemma \ref{lemma:adr-input}. Then, the spaces $X$  and
  $X_{\nb}$ are replaced by $\tX$ and $\tX_{\tnb}$. The basis $\{\tphi_1, \dots,
  \tphi_{\tnb}\}$ is equal to $\{\phi_1, \dots, \phi_{\nb}\}$ with the addition of a constant
  function, which can be emulated exactly by deep ReLU neural networks. 
  Hence, Proposition \ref{prop:basis} can be extended to this case.
  Finally, the matrices $\AaNbNq$ and $\AoneNbNq$ used in the proof of Theorem
  \ref{th:deepONet} are replaced, respectively, by the matrices with entries
  \begin{equation*}
    \frab_{\nq}^{(\argsADR)}(\tphi_j, \tphi_i) \quad\text{and}\quad \frab_{\nq}^{(\Id_d, 1)}(\tphi_j, \tphi_i), \qquad (i, j)\in \brange{\tnb}^2,
  \end{equation*}
  where
  \begin{equation*}
    \frab_{\nq}^{(\argsADR)} (u, v)\coloneqq \sum_{k=1}^{\nq } \wkq \left(\bA(\xkq) \nabla u(\xkq)\right) \cdot \nabla v (\xkq) + \sum_{k=1}^{\nq}\wkq c(\xkq) u(\xkq) v(\xkq), 
  \end{equation*}
  for all $u, v\in C^1(\Omega)$.
 Since the bilinear form $\bADR$ is coercive and continuous on $H^1(\Omega)$,
 results equivalent to Lemmas \ref{lemma:contcoer} and \ref{lemma:norminvA} with
 the new matrices can be proven directly. The rest of the proof is the same as
 the proof of Theorem \ref{th:deepONet}.
\end{proof}
\subsection{Linear Elastostatics}
\label{sec:Elast}
%
\subsubsection{Definition of the problem}
\label{sec:Elast-prob}
We assume $d=2,3$. 
Small, linear elastic deformation of a body occupying $Q =(0,1)^d$ 
with periodic boundary conditions and
subject to a prescribed, periodic body force $\bff:\Omega = \R^d/\Z^d \to \mathbb{R}^d$
can be described by the displacement field $\bu:\Omega \to \mathbb{R}^d$ 
which satisfies the equilibrium of stress
\begin{equation}\label{eq:LinElPDE}
\begin{aligned}
\div\bsigma[\bu] + \bff &=0 \quad\text{in }\Omega.
\end{aligned}
\end{equation}
Here $\bsigma[\bu]:\Omega \to \Rsymdd$ is a symmetric matrix function,
the \emph{stress tensor}. 
It depends on the displacement field $\bu$
via the (linearized) \emph{strain tensor}
$\bepsilon[\bu]:\Omega\to \Rsymdd$ which is given by
\begin{equation}\label{eq:LinElStrain}
\bepsilon[\bu] \coloneqq \frac{1}{2}\left(\grad \bu + (\grad \bu)^\top\right),
\quad
(\bepsilon[\bu])_{ij} \coloneqq \frac{1}{2}(\partial_j u_i + \partial_i u_j),\quad i,j=1,...,d\;.
\end{equation}
In the linearized theory,
the tensors $\bsigma$ and $\bepsilon$ in \eqref{eq:LinElPDE}, \eqref{eq:LinElStrain}
are related by the linear constitutive stress-strain relation (``Hooke's law'')
\begin{equation}\label{eq:Hooke}
\bsigma = \ttA \bepsilon .
\end{equation}
In \eqref{eq:Hooke}, $\ttA$ is a fourth order tensor field,
i.e. $\ttA = \{ \ttA_{ijkl} : i,j,k,l=1,...,d\}$,
with certain symmetries: the $d^4$ component functions $\ttA_{ijkl}(x)$
are assumed analytic in $[0,1]^d$ and $1$-periodic with respect to each coordinate,
and satisfy for every $x\in \Omega$,
\begin{equation}\label{eq:Asym}
\forall \btau \in \Rsymdd, \;\ttA(x) \btau \in \Rsymdd
\quad \text{and}\quad
\forall \btau,\bsigma \in \Rsymdd, \; (\ttA(x) \btau ): \bsigma = (\ttA(x) \bsigma) : \btau \;.
\end{equation}
Key assumption on $\ttA$ is \emph{coercivity}:
there exists a constant $\amin>0$ such that
\begin{equation}\label{eq:Acoerc}
\forall x\in \Omega, \; \forall \btau \in \Rsymdd, 
\quad 
(\ttA(x) \btau) : \btau \geq \amin \| \btau \|^2_2  \;.
\end{equation}
see, e.g., \cite{Truesdell} for details. 
Inserting \eqref{eq:Hooke} into \eqref{eq:LinElPDE},
integrating by parts and noting the periodic boundary conditions,
the ``primal variational formulation'' of \eqref{eq:LinElPDE} 
reads: 
find $\bu^{\ttA}\in [H^1(\Omega)/\mathbb{R}]^d$ such that 
\begin{equation}\label{eq:ElasBilin}
\frab^\ttA(\bu^{\ttA},\bv) 
\coloneqq \int_\Omega \eps[\bv]: (\ttA \eps[\bu^{\ttA}]) 
= \int_\Omega \bbf\cdot \bv \qquad \forall \bv\in [H^1(\Omega)/\mathbb{R}]^d\;.
\end{equation}
Unique solvability of \eqref{eq:ElasBilin} is implied by the Lax-Milgram Lemma
with \eqref{eq:Acoerc} and Korn's inequality 
upon noticing that the space $X^d = [H^1(\Omega)/\mathbb{R}]^d$ 
does not contain rigid body motions: 
rigid body rotations are eliminated due to the periodicity of the present setting,
and 
rigid body translations with the factoring of constants in each component.
The Korn inequality and the Poincaré inequality \eqref{eq:poincare} imply
existence of a positive constant $c$ such that
\begin{equation*}
\forall \bv \in X^d: \quad \frab^\ttA(\bv,\bv) \geq c \amin \| \bv \|_{H^1(Q)}^2 \;.
\end{equation*}
For given, fixed, $Q$-periodic
$\bbf\in [\Hol(\Omega)/\mathbb{R}]^d$, 
there exists a unique solution of \eqref{eq:ElasBilin}. Furthermore,
the coefficient-to-solution map $\Sel : \ttA \mapsto \bu^{\ttA}$ is analytic from the set  
$\cDel = \{ \ttA \in \Hol(\Omega, A_{\cDel})^{d^4}: \text{\eqref{eq:Acoerc} and \eqref{eq:Asym} hold} \}$
to 
$\cUel = \Sel(\cDel) \subset X^d\cap \Hol(\Omega, A_{\cUel})^d$, for positive constants $A_{\cDel}, A_{\cUel}$.
\subsubsection{Operator network approximation}
\label{sec:ONetApr}
For the operator network approximation of the map $\Sel$, we introduce modified
encoding and reconstruction operators. To construct the encoding operator, we
extend the definition of the vectorization operation to fourth order tensors so
that, for all $\ttB\in R^{n_1\times \dots \times n_4}$, $\vec(\ttB)\in
\R^{n_1\cdots n_4}$.
We consequently extend the definition of the reordering function introduced in
Section \ref{sec:ADR} to $\frav : \N\times\N\times\N\times\N$ such that
\begin{equation}
  \label{eq:frav-2}
  \vec(\ttB)_{\frav(m,n,p,q)} = \ttB_{mnpq}.
\end{equation}
The modified encoding operator $\cEel : [C([0,1]^d)]^{d\times d\times d\times d}
\to \R^{d^4\nq}$ is then given by
\begin{equation}
  \label{eq:Eel}
  \cEel(\ttA) \coloneqq
  \begin{pmatrix}
    \vec(\ttA(\xq_1))\\
    \vdots \\
    \vec(\ttA(\xq_{\nq}))
  \end{pmatrix},
\end{equation}
where $\xenc = \xq_1, \dots, \xq_{\nq}$ are the usual quadrature points.
For all $m\in \N$, the modified reconstruction operator $\cRel : \R^{dm} \to
H^1(\Omega)^d$ is instead defined, 
given a neural network $\Phibranchnoeps$ such that $\Realiz(\Phibranchnoeps) : \overline{Q}\to
\R^m$, as
\begin{equation}
  \label{eq:Rel}
  \cRel_{\Phibranchnoeps}(\bc)(x) = \left( \Id_d \otimes \Realiz(\Phibranchnoeps)(x) \right)^\top \bc, \qquad \forall x\in \overline{Q},\; \forall \bc\in \R^{dm}.
\end{equation}
We can now state the operator network approximation result for problem \eqref{eq:LinElPDE}.
\begin{theorem}
  \label{th:el}
  Theorem \ref{th:deepONet} holds with $a\in\cD$ replaced by $\ttA\in\cDel$,
  $S(a)$ replaced by $\Sel(\ttA)$, $\cE_{\xenc}(a)$ replaced by $\cEel(\ttA)$, and $\cR_{\Phibranch}$
  replaced by $\cRel_{\Phibranch}$.
\end{theorem}
\begin{proof}
  We construct a basis of the $d\nb$-dimensional discrete space $X_{\nb}^d$ approximating
  $X^d$ as
  \begin{equation*}
    \bpsi_1 =
    \begin{pmatrix}
      \phi_1 \\
      \vdots\\
      0
    \end{pmatrix},
    \dots,
\bpsi_{\nb} =
    \begin{pmatrix}
      \phi_{\nb} \\
      \vdots\\
      0
    \end{pmatrix},
\dots,
\bpsi_{(d-1)\nb+1} =
    \begin{pmatrix}
      0 \\
      \vdots\\
      \phi_{1}
    \end{pmatrix},
\bpsi_{d\nb} =
    \begin{pmatrix}
      0 \\
      \vdots\\
      \phi_{\nb}
    \end{pmatrix},
  \end{equation*}
  where $\phi_1, \dots, \phi_{\nb}$ are the basis functions defined in
  Section \ref{sec:pol-quad}.
  The trunk network $\Phitrunk$ is then constructed as in the proof of Theorem
  \ref{th:deepONet}: it follows that the $j$th column of
  \[
    \left(\Id_d \otimes \Realiz(\Phitrunk) \right)^\top
  \]
  contains an approximation of $\bpsi_j$, for each $j\in \brange{d\nb}$.

  To construct the branch network $\Phibranch$, we replace the input layer
  used in the proof of Theorem \ref{th:deepONet}, in a similar way as we did in Lemma
  \ref{lemma:adr-input}. Define,
  for all $i,j\in \brange{d\nb}$ and $m,n,p,q\in \brange{d}$,
  \begin{equation*}
    \tbD^{ij}_{mnpq}(\xkq) = \wkq\left(\bepsilon[\bpsi_i](\xkq)  \right)_{mn}\left(\bepsilon[\bpsi_j] (\xkq) \right)_{pq}
  \end{equation*}
  and let $\tvec(\tbD(\xkq))\in \R^{d^2\nb^2\times d^4}$ such that
  \begin{equation*}
    \tvec(\tbD(\xkq))_{\frav(i,j)\frav(m,n,p,q)} = \tbD^{ij}_{mnpq}(\xkq),
  \end{equation*}
  with $\frav$ defined in \eqref{eq:frav-1} and \eqref{eq:frav-2} for two and
  four arguments, respectively. Then,
  \begin{equation*}
    \Phiinputel
  \coloneqq \left(
    \left(-\alpha \left[\tvec(\tbD(\xq_1)) | \dots |\tvec(\tbD(\xq_{\nq }))\right], \bzero_{d^2\nb^2}\right)
  \right)
  \end{equation*}
  is such that
  \begin{equation*}
    \matr\left(\Realiz(\Phiinputel) (\cEel(\ttA))   \right)_{ij} = 
  -\alpha \frab^{\ttA} (\bpsi_j, \bpsi_i), \qquad \forall (i, j)\in \brange{d\nb}^2.
  \end{equation*}
  We can then construct $\Phitrunk$ as in the proof of Theorem
  \ref{th:deepONet}, with $\Phiinputel$ replacing $\Phiinputone$. The rest of
  the proof follows the same argument as the proof of Theorem \ref{th:deepONet}.
\end{proof}
\section{Conclusions} 
\label{sec:conc}
We proved, in the periodic setting on $\Omega = \R^d/\Z^d $, 
the exponential convergence of deep operator network emulation 
of the coefficient-to-solution map of some linear elliptic equations, 
under the assumption of analytic coefficients $a$ and right-hand sides $f$.
The proof used the analytic regularity of solutions 
$u^a$ of \eqref{eq:problem}
implied by classical elliptic regularity results and
the consequential exponential convergence of
polynomial approximations of $a$ and $u^a$ and of
fully discrete spectral-Galerkin numerical schemes.
The expression rate bounds were not explicit in the physical domain dimension $d$
which is moderate in engineering applications.
We have developed the analysis for isotropic diffusion equations 
and
extended it to problems with parametric diffusion, with anisotropic diffusion
and reaction, and to linear elastostatics.
We also obtained corresponding expression rate bounds for PDEs with parametric
inputs. Here, leveraging NN composition, dependence of NN expression rates 
for the parametric inputs on the parameter dimension $\dimp$ is inherited by
the solution expression rate bounds.

A direct application of the present results are expression rate bounds
for homogenized (``upscaled'') coefficients for linear, elliptic 
divergence-form PDEs with diffusion coefficients of the form $a_{ij}(x/\eta)$
where $x\in \Omega$, $0<\eta \ll {\rm diam}(\Omega)1$ 
is a small length scale, and 
$y\mapsto a_{ij}(y)$ is analytic and $Y$-periodic in the unit cell $Y=[0,1]^d$.
As it is well-known (e.g. \cite[Sec.~6]{CDHomog99}), 
in space dimension $d\geq 2$, solutions $u^\eta\in H^1_0(\Omega)$ of, e.g., 
the homogeneous Dirichlet problem in a bounded domain $\Omega \subset \R^d$ of
unit size converge, as $\eta \downarrow 0$, 
in $L^2(\Omega)$ to the solution $u^0$
of a linear, elliptic PDE with ``upscaled'' coefficients $a_{ij}^0$ that are
constant in $\Omega$.
These coefficients admit simple expressions (see, e.g., \cite[Prop.~6.8]{CDHomog99}) 
in terms of solutions $w_j(y)$ of a linear, elliptic, $Y$-periodic boundary value
problem with coefficients and right-hand sides analytic in $Y$.
The ``unit-cell solution map'' $S:(a_{ij}(y): 1\leq i,j \leq d) \mapsto (w_j(y):j=1,...,d)$ 
is covered by the results of the present paper and can, hence, be 
approximately expressed with accuracy 
$\varepsilon>0$ 
by an ONet of size polylogarithmic w.r.t.~$\varepsilon$.
The branch part of the ONet approximating $S$ can also be concatenated with an
averaging network (performing the integration of \cite[Prop.~6.8]{CDHomog99})
to construct an ``upscaling
NN'' emulating the map
$\widehat{S}: (a_{ij}(y): 1\leq i,j \leq d) \mapsto (a_{ij}^0:  1\leq i,j \leq
d)$.
With the results of Section~\ref{sec:Elast},
corresponding statements hold for elasticity based on 
\cite[Chap.~10]{CDHomog99} and \cite[Prop.~10.12]{CDHomog99}.

Further directions comprise the derivation of NN expression rate bounds 
of data-to-solution maps for \emph{smooth, nonlinear forward problems} 
and for \emph{inverse problems}.
In the analysis of such problems, the presently established NN expression rate
bounds constitute an essential building block in the construction of parameter-sparse
NN approximations.
Our analysis extends to ONet approximations of operators 
    mapping both coefficients and right-hand sides to the solutions, see Remark \ref{remark:fmap}.
\appendix
\section{Exponential convergence of fully discrete Spectral-Galerkin Solution}
\label{sec:GNI}
We present here the exponential convergence of fully discrete
  Spectral-Galerkin solutions of the problems considered in this paper. 
The following classical approximation result will be useful.
\begin{lemma}
  \label{lemma:analytic-approx}
  Let $A>0$. Let $X= \{v\in H^1(\Omega): \int_Q v = 0\}$ or
    $X=H^1(\Omega)$. 
  Then, there exist constants $C, b>0 $ such that for all $p\in\N_0$ and 
  for all $v\in\Hol(\Omega; A)$,
  \begin{equation*}
   \inf_{w\in \Q_p(Q)\cap X} \| v-w\|_{H^1(Q)}  \leq C\exp(-bp),
   \qquad
   \inf_{w\in \Q_p(Q)\cap X} \| v-w\|_{L^\infty(Q)}  \leq C\exp(-bp).
  \end{equation*}
\end{lemma}
\begin{proof}
  The error bound in the $H^1(Q)$-norm follows from tensorization of a
    univariate interpolation operator that is exact at the endpoints, 
    see, e.g., \cite[Theorem A.2]{Costabel2005}. This implies
    \begin{equation*}
      \inf_{w\in \Q_p(Q)\cap H^1(\Omega)} \| v-w\|_{H^1(Q)}  \leq C_1\exp(-bp).
    \end{equation*}
  The bound in the $L^\infty(Q)$-norm is a consequence of, e.g., 
  \cite[Remark 3.1 and Theorem 3.5]{Opschoor2019}, which implies
  \begin{equation*}
    \inf_{w\in \Q_p(Q)} \| v-w\|_{L^\infty(Q)}  \leq C_1\exp(-bp).
  \end{equation*}
  To ensure conformity in $H^1(\Omega) = \Hper^1(Q)$, i.e., to impose continuity
  across matching hyperfaces, it is sufficient to lift the difference at
  vertices, followed (if $d\geq 2$) by edges, and iteratively up to $d-1$
  dimensional hyperfaces. The norm of the lifting is bounded by a constant
  (exponential in $d$) multiplying the $L^\infty(Q)$-norm of $v-w$.

  If $X$ includes the vanishing average constraint, it is sufficient to remark that
  for all $w \in L^2(Q)$ and all $v\in X$, 
  it holds that $\left| \int_Q w  \right|\leq \| v-w\|_{L^2(Q)}$.
\end{proof}
  The following lemma, then, concerns the convergence of fully discrete Spectral-Galerkin
  solutions for problems in $\Omega$, with analytic right-hand sides and coefficients. 
\begin{lemma}
 \label{lemma:GNI} 
 Let $\genset \in \{\cD, \tcD, \cDel\}$ and $\frad=d$ for linear elasticity,
 $\frad=1$ otherwise.
 Let $X$ be the space of solutions of the problem considered and denote
   $X_p = \Q_p(Q)^\frad\cap X$.
 Let $f\in \Hol(\Omega)^{\frad}$
 and, for coefficients
 $\genpar\in\genset$, let
 $\frab^{\genpar}(\cdot, \cdot)$
   be one of the bilinear forms defined in Sections
   \ref{sec:problem}, \ref{sec:ADR-prob}, or \ref{sec:Elast-prob}.
 There exists $C_1, C_2>0$ such that, for all $p\in \N$
 and for all integer $q\geq p+1$, 
 \begin{equation*}
   \sup_{\genpar\in\genset}\|u^{ \genpar } - \ugenparNbNq\|_{H^1(Q)} \leq C_1 \exp (-C_2p),
 \end{equation*}
 where $\ugenparNbNq\in X_p$ 
 is such that 
 $\frab^{\genpar}_{q^d}(\ugenparNbNq, v)= (f, v)$ for all $v\in X_p$.
\end{lemma}
\begin{proof}
Strang's lemma \cite[Lemma 10.1]{Quarteroni2017}
implies that
there exists $C>0$ independent of $\genpar\in \genset$, $p$, 
and $q$, such that
\begin{equation*}
    \| u^{\genpar} - \ugenparNbNq \|_{H^1(Q)} \leq C \inf_{v\in X_p} \left(\| u^{\genpar} - v \|_{H^1(Q)}
      + \sup_{w\in X_p  \setminus \{\bzero\}} \frac{|\frab^{\genpar}(v, w) - \frab^{\genpar}_{\nq}(v, w)|}{\|v\|_{H^1(Q)}\|w\|_{H^1(Q)}}
  \right).
\end{equation*}
By \cite[Section 6.4.3]{CHQZ1}, then, denoting $\tp = \lfloor p/2 \rfloor$,
there exists $\tC>0$ independent of $\genpar\in\genset$, $p$, and $q$ such that
\begin{equation*}
    \| u^{\genpar} - \ugenparNbNq \|_{H^1(Q)} 
    \leq \tC  \left(\inf_{v\in X_p }\| u^{\genpar} - v \|_{H^1(Q)}
      +  \inf_{v\in Y_{\tp}} \| \genpar- v\|_{L^\infty(Q)} \right),
\end{equation*}
where the space $Y_{\tp}$ depends on the problem under consideration:
\begin{equation*}
    Y_{\tp} =
    \begin{cases}
     X_{\tp}  &\text{if } \genset = \cD,\\
   \left\{\bA \in X_{\tp}^{d\times d}: 
        \bA_{ij} = \bA_{ji}\right\} \times X_{\tp} &\text{if } \genset = \tcD,
    \\
    \left\{\ttA \in X_{\tp}^{d\times d\times d\times d}: 
    \text{\eqref{eq:Asym} holds}\right\} &\text{if } \genset = \cDel.
    \end{cases}
  \end{equation*}
  Since functions in $\genset$ and in $S(\cD)$, $\SADR(\tcD)$, or $\Sel(\cDel)$ are analytic with uniform bounds
  on the norms at all orders, using Lemma \ref{lemma:analytic-approx} concludes the proof.
 \end{proof}
\section{Lipschitz continuity of the data-to-solution map}
\label{sec:LipSolOp}
For the readers' convenience, 
    we provide a proof of the (known) Lipschitz dependence of the solution of
    the PDEs considered in this paper
    on the coefficients.
\begin{lemma}
  \label{lemma:lip-sol}
  Let $X$ be a Hilbert space, let $Y$ be a Banach space, and let
  $\genset\subset Y$. 
  Let furthermore $\frab^{\genpar}: X\times X\to \R$  be a bilinear
  form that is also linear with respect to the coefficient $\genpar$. 
  Suppose that there exists $C_{\mathrm{cont}}>0$ such that
  \begin{equation}
    \label{eq:gencont}
 \forall \genpar\in Y: \qquad 
    \frab^{\genpar} (u, v)  \leq
    C_{\mathrm{cont}}
    \| \genpar \|_Y \| u \|_X \| v\|_X, \qquad \forall u, v\in X\;.
  \end{equation}
  Furthermore, suppose there exists $\genparmin>0$ such that
  \begin{equation}
    \label{eq:gencoer}
    \frab^{\genpar}(u, u )\geq \genparmin \|u\|_X^2,\qquad \forall u\in X,\;\forall \genpar\in \genset.
  \end{equation}
  For fixed $f\in X'$ and for each $\genpar\in \genset$, define $u^{\genpar}\in X$ as
  the function such that
  \begin{equation}
    \label{eq:gensol}
    \frab^{\genpar}(u^{\genpar}, v) = \langle f, v \rangle, \qquad \forall v\in X.
  \end{equation}
  Then, there exists $C_L>0$ (depending only on $C_{\mathrm{cont}}$, $\genparmin$, and $f$) 
such that
  \begin{equation*}
    \| u^{\genpar_1} - u^{\genpar_2} \|_{X} \leq C_L \| \genpar_1 -\genpar_2\|_{Y}, \qquad \forall \genpar_1, \genpar_2\in \genset.
  \end{equation*}
\end{lemma}
\begin{proof}
Denote $u_i = u^{\genpar_i}$, $i=1,2$. Using \eqref{eq:gencoer}, \eqref{eq:gensol}, the continuity of the bilinear form
with respect to the coefficient, and \eqref{eq:gencont}
  \begin{align*}
    \| u_1 - u_2 \|_X^2
    &\leq \frac{1}{\genparmin} \frab^{\genpar_1}(u_1-u_2, u_1-u_2)
    \\ &\leq
         \frac{1}{\genparmin}\left(\frab^{\genpar_2}(u_2, u_1-u_2) - \frab^{\genpar_1}(u_2, u_1-u_2)  \right)
    \\ &\leq \frac{C_{\mathrm{cont}}}{\genparmin}\| \genpar_2 - \genpar_1 \|_{Y} \| u_2\|_{X} \| u_1 - u_2\|_{X}.
  \end{align*}
  The Lax-Milgram bound $ \| u_2\|_{X} \leq \frac{1}{\genparmin} \|f\|_{X'} $
  concludes the proof.
\end{proof}
\section{Generalized eigenvalue problems}
We recall and a result on the relationship between Rayleigh
quotients and generalized eigenvalue problems (see, e.g., \cite[Section I.10]{Strang2019}).
\begin{lemma}
  \label{lemma:eigprob}
  Let $n\in\N$ and let $\bA, \bB \in \R^{n\times n}$ be symmetric, positive definite
  matrices. 
  Suppose that there exist constants $c, C>0$ such that
  \begin{equation*}
    \forall \bx\in \R^n: \qquad 
    c \leq \frac{\bx^\top \bA \bx}{\bx^\top \bB \bx} \leq C.
  \end{equation*}
  Then, the spectrum of $\rev{\bB^{-1/2}\bA\bB^{-1/2}}$ is contained in $[c, C]$.
\end{lemma}
\begin{proof}
  Denote by $\bB^{1/2}$ the symmetric positive definite matrix such that
  $\bB^{1/2}\bB^{1/2}=\bB$ and by $\bB^{-1/2}$ its inverse. For all $\by\in\R^n$,
  we can write $\bx = \bB^{-1/2}\by$ and by hypothesis
  \begin{equation*}
    \frac{\by^\top\bB^{-1/2}\bA\bB^{-1/2} \by}{\by^\top\by} = 
    \frac{\bx^\top \bA \bx}{\bx^\top \bB \bx}  \in [c, C].
  \end{equation*}
  Now, $\bB^{-1/2}\bA\bB^{-1/2}$ is a symmetric matrix, and we have shown that
  its Rayleigh quotient is contained in $[c, C]$. Therefore,
    $\sigma(\bB^{-1/2}\bA\bB^{-1/2}) \subset [c, C]$.
\end{proof}
\section{Polynomial basis in $\Hper^1(Q)$}
  \label{sec:L2l2}
  Let $X_{\nb}$ be defined as in \eqref{eq:Xnb}; we adopt the notation of
  Section \ref{sec:pol-quad}. In the next lemma, we show that $X_{\nb}$ is a
  basis for polynomials that is conforming in $H^1(\Omega)$.
  \begin{lemma}
    \label{lemma:basis}
   Let $X = \{v\in H^1(\Omega) = \Hper^1(Q) :  \int_Q v = 0\}$. Let $p\in \N$
   and $\nb = p^d-1$. Then,
   \begin{equation*}
     X_{\nb} = \Q_p(Q)\cap X.
   \end{equation*}
  \end{lemma}
  \begin{proof}
    The inclusion $X_{\nb}\subset \Q_p(Q)\cap X$ is a consequence of the
    definition of $\phi_1, \dots, \phi_{\nb}$. We now prove that the dimensions
    of $X_{\nb}$ and $\Q_p(Q)\cap X$ are equal.
     Set $\zetaoned_0 = L_0$,
    $\zetaoned_1 = L_1$, and
    \begin{align*}
      \zetaoned_{2i} &= L_{2i}, \qquad i =1, \dots, \lfloor p/2 \rfloor,\\
      \zetaoned_{2i+1}&= L_{2i+1} - L_1, \qquad i = 1, \dots, \lfloor (p-1)/2 \rfloor.
    \end{align*}
    Since $\{L_0, \dots, L_p\}$ is a basis for $\Q_p((0,1))$, then
    $\{\zetaoned_0, \dots, \zetaoned_p\}$ is also a basis for $\Q_p((0,1))$.
    Define
    \begin{equation*}
      \zeta_\is = \zetaoned_{i_1}\otimes \dots \otimes \zetaoned_{i_d}, \qquad (\is)\in \N_0^d.
    \end{equation*}
    It follows that $\{\zeta_{\is}: (\is)\in \{0, \dots, p\}^d\}$ is a basis for $\Q_p(Q)$.
    Now, given any function
    \begin{equation*}
      v = \sum_{\is\in \{0, \dots, p\}^d} v_\is \zeta_\is,
    \end{equation*}
    requiring that $v\in \Hper^1(Q)$ is equivalent to imposing that
    \begin{equation*}
      v_\is = 0, \qquad \forall ( \is ) \in \{0, \dots, p\}^d: \exists j: i_j=1.
    \end{equation*}
    Since $\dim(\Q_p(Q)) = (p+1)^d$ and
    \begin{equation*}
      \card\left(\left\{ ( \is ) \in \{0, \dots, p\}^d: \exists j: i_j=1 \right\}  \right) = \sum_{n=1}^d \binom{d}{n}p^{d-n}, 
    \end{equation*}
    we infer that 
    \begin{equation*}
      \dim(\Q_p(Q) \cap \Hper^1(Q)) = (p+1)^d - \sum_{n=1}^d \binom{d}{n}p^{d-n} = p^d.
    \end{equation*}
    Imposing vanishing average removes another degree of freedom, which
    implies $\dim(\Q_p(Q)\cap X) = \nb$ and concludes the proof.
  \end{proof}
  In the next lemma we bound the $\ell^2$ norm of the coefficients (in the
  basis introduced in Section \ref{sec:pol-quad}) of a
  function in $X_{\nb}$ by its $L^2$-norm.

 \begin{lemma}
    \label{lemma:L2l2}
    Let $\phi_i$, $i \in \N_0$, be the functions defined in \eqref{eq:phibasis}.
   There exists $\CLell \geq 1$ such
   that, for all $p\in \N$ with $p\geq 2$ and for all $v\in X_{\nb}$ with  $\nb =
   p^d-1$ such that $v = \sum_{i=1}^{\nb} w_i \phi_i$,
  \begin{equation}
    \label{eq:L2-l2}
    \| \bw \|_2\leq \CLell \nb^{1/2}\| v \|_{L^2(Q)} 
  \end{equation}
  where  $\bw = (w_1, \dots, w_{\nb})$.
\end{lemma}
\begin{proof}
  For ease of notation, we introduce the functions
  \begin{equation*}
   \nuoned_{2i} = L_{2i}, \qquad \nuoned_{2i+1} = L_{2i+1}-L_1, \qquad i\in \N_0,
  \end{equation*}
  and
  \begin{equation*}
    \nu_{i_1, \dots, i_d} = \nuoned_{i_1}\otimes \dots\otimes \nuoned_{i_d}, \qquad (i_1, \dots, i_d)\in\N_0^d.
  \end{equation*}
  Remark that $\nuoned_1= 0$. We then rewrite $v$ as
  \begin{equation*}
    v = \sum_{(i_1, \dots, i_d)\in \{0, \dots, p\}^d} x_{i_1, \dots, i_d} \nu_{i_1, \dots, i_d}.
  \end{equation*}
  Note that, since $v\in X_{\nb}$, we have $x_{0, \dots, 0} =0$. Furthermore, we
  assume that $x_{i_1, \dots, i_d} = 0$ if there exists at least one $j\in \{1,
  \dots, d\}$ such that $i_j=1$.
  We have 
  \begin{equation}
    \label{eq:x-w}
    \left\{x_{i_1, \dots, i_d} : (i_1, \dots, i_d)\neq (0, \dots, 0) \text{ and }\forall j\in \{1, \dots, d\}, i_j\neq 1\right\} = \{w_1, \dots, w_{\nb}\},
  \end{equation}
  with one-to-one correspondence between the sets.
  We now introduce the $L^2((0,1))$-orthonormal basis
  \begin{equation*}
    \chioned_i = \frac{L_i}{\|L_i\|_{L^2((0,1))}}, \qquad i\in \N_0,
  \end{equation*}
  and the $L^2(Q)$-orthonormal basis
  \begin{equation*}
    \chi_{i_1, \dots, i_d} = \chioned_{i_1}\otimes \dots \otimes \chioned_{i_d}, \qquad (i_1, \dots, i_d) \in \N_0^d.
  \end{equation*}
  It follows that
  \begin{equation}
    \label{eq:L2chi}
    \begin{aligned}
      \| v\|_{L^2(Q)}^2
      &=
        \sum_{(i_1, \dots, i_d)\in \{0, \dots, p\}^d} (v, \chi_{i_1, \dots, i_d})_{L^2(Q)}^2 \\
      & = 
        \sum_{(i_1, \dots, i_d)\in \{0, \dots, p\}^d} \left(\sum_{(j_1, \dots, j_d)\in \{0, \dots, p\}^d}x_{j_1, \dots, j_d} \prod_{k=1}^d(\nuoned_{j_k}, \chioned_{i_k})_{L^2((0,1))}  \right)^2.
    \end{aligned}
  \end{equation}
  Now, for $i, j\in \{0, \dots, p\}$, we have
  \begin{equation*}
    (\nuoned_j, \chioned_i)_{L^2((0,1))} =
    \begin{cases}
      \delta_{ij} \| L_i \|_{L^2((0,1))} &\text{if } i\neq 1\\
      - \| L_1 \|_{L^2((0,1))} &\text{if }i=1\text{ and $j$ is odd, $j\neq 1$}\\
      0 &\text{otherwise},
    \end{cases}
  \end{equation*}
  where $\delta_{ij}$ is the Kronecker delta.
  Hence,
  \begin{multline*}
    \|v\|^2_{L^2(Q)} 
    = \sum_{(i_1, \dots i_d): i_j\neq 1, \forall j} x_{i_1, \dots, i_d}^2\prod_{k=1}^d \|L_{i_k}\|^2_{L^2((0, 1))}
   \\ +
    \sum_{(i_1, \dots i_d): \exists j \text{ s.t. }i_j= 1}
    \left(\sum_{(j_1, \dots, j_d)\in \{0, \dots, p\}^d}x_{j_1, \dots, j_d} \prod_{k=1}^d(\nuoned_{j_k}, \chioned_{i_k})_{L^2((0,1))}  \right)^2,
  \end{multline*}
  where in the sums we still require also $(i_1, \dots, i_d)\in \{0, \dots, p\}^d$.
  Since $\|L_i \|^2_{L^2((0,1))} = 1/(2i+1)$, it follows from the equation above that
  \begin{equation*}
    \|v\|^2_{L^2(Q)} \geq \frac{1}{(2p+1)^{d}} \sum_{(i_1, \dots i_d): i_j\neq 1, \forall j} x_{i_1, \dots, i_d}^2.
  \end{equation*}
 From the correspondence between coefficients \eqref{eq:x-w} and since $\nb =
 p^d-1$, we obtain the bound in \eqref{eq:L2-l2}. We remark that the constant
 $\CLell$ depends exponentially on $d$.
\end{proof}
\bibliographystyle{siamplain}
\bibliography{library}

@book{CDHomog99,
    AUTHOR = {Cioranescu, Doina and Donato, Patrizia},
     TITLE = {An introduction to homogenization},
    SERIES = {Oxford Lecture Series in Mathematics and its Applications},
    VOLUME = {17},
 PUBLISHER = {The Clarendon Press, Oxford University Press, New York},
      YEAR = {1999},
     PAGES = {x+262},
      ISBN = {0-19-856554-2},
   MRCLASS = {35B27 (35J20 35K05 35L05 49J45 74-02 74Qxx 78M40)},
  MRNUMBER = {1765047},
MRREVIEWER = {Eugenia P\'{e}rez},
}

@article {Costabel2005,
    AUTHOR = {Costabel, Martin and Dauge, Monique and Schwab, Christoph},
     TITLE = {Exponential convergence of {$hp$}-{FEM} for {M}axwell
              equations with weighted regularization in polygonal domains},
   JOURNAL = {Math. Models Methods Appl. Sci.},
  FJOURNAL = {Mathematical Models and Methods in Applied Sciences},
    VOLUME = {15},
      YEAR = {2005},
    NUMBER = {4},
     PAGES = {575--622},
      ISSN = {0218-2025},
   MRCLASS = {78M10 (65N30)},
  MRNUMBER = {2137526},
MRREVIEWER = {Snorre H. Christiansen},
       DOI = {10.1142/S0218202505000480},
       URL = {https://doi.org/10.1142/S0218202505000480},
}

@ARTICLE{DeRyck2022,
       author = {{De Ryck}, Tim and {Mishra}, Siddhartha},
        title = "{Generic bounds on the approximation error for physics-informed (and) operator learning}",
      journal = {arXiv e-prints},
         year = 2022,
        month = may,
          eid = {arXiv:2205.11393},
        pages = {arXiv:2205.11393},
archivePrefix = {arXiv},
       eprint = {2205.11393},
 primaryClass = {cs.LG},
       adsurl = {https://ui.adsabs.harvard.edu/abs/2022arXiv220511393D}
}

@ARTICLE{Li2020,
       author = {{Li}, Zongyi and {Kovachki}, Nikola and {Azizzadenesheli}, Kamyar and {Liu}, Burigede and {Bhattacharya}, Kaushik and {Stuart}, Andrew and {Anandkumar}, Anima},
        title = "{Fourier Neural Operator for Parametric Partial Differential Equations}",
      journal = {arXiv e-prints},
         year = 2020,
        month = oct,
          eid = {arXiv:2010.08895},
        pages = {arXiv:2010.08895},
archivePrefix = {arXiv},
       eprint = {2010.08895},
 primaryClass = {cs.LG},
       adsurl = {https://ui.adsabs.harvard.edu/abs/2020arXiv201008895L}
}

@article{Kovachki2021,
  author  = {Nikola Kovachki and Samuel Lanthaler and Siddhartha Mishra},
  title   = {On Universal Approximation and Error Bounds for Fourier Neural Operators},
  journal = {Journal of Machine Learning Research},
  year    = {2021},
  volume  = {22},
  number  = {290},
  pages   = {1--76},
  url     = {http://jmlr.org/papers/v22/21-0806.html}
}

@article{deng2021convergence,
title = {{Approximation rates of DeepONets for learning operators arising from advection-diffusion equations}},
journal = {Neural Networks},
year = {2022},
issn = {0893-6080},
doi = {https://doi.org/10.1016/j.neunet.2022.06.019},
author = {Beichuan Deng and Yeonjong Shin and Lu Lu and Zhongqiang Zhang and George Em Karniadakis},
}

@book{strang2019,
  title={Linear algebra and learning from data},
  author={Strang, Gilbert},
  volume={4},
  year={2019},
  publisher={Wellesley-Cambridge Press Cambridge}
}

@article{Grohs2021,
    author = {Grohs, Philipp and Herrmann, Lukas},
    title = "{Deep neural network approximation for high-dimensional elliptic {PDEs} with boundary conditions}",
    journal = {IMA Journal of Numerical Analysis},
    year = {2021},
    month = {05},
    issn = {0272-4979},
    doi = {10.1093/imanum/drab031},
}

@Techreport{HSZ,
    author = {Lukas Herrmann and Christoph Schwab and Jakob Zech},
    title = {Neural and gpc operator surrogates: construction and expression rate bounds},
    institution = {Seminar for Applied Mathematics, ETH Z{\"u}rich},
    number = {2022-27},
    year = {2022}
}

@book{Truesdell,
    AUTHOR = {Truesdell, III, Clifford A.},
     TITLE = {A first course in rational continuum mechanics. {V}ol. 1},
    SERIES = {Pure and Applied Mathematics},
    VOLUME = {71},
   EDITION = {Second},
      NOTE = {General concepts},
 PUBLISHER = {Academic Press, Inc., Boston, MA},
      YEAR = {1991},
     PAGES = {xviii+391},
      ISBN = {0-12-701300-8},
   MRCLASS = {73A05 (73-01 73B05 76A02)},
  MRNUMBER = {1162744},
MRREVIEWER = {G\'{e}rard A. Maugin},
}

@book {CHQZ1,
    AUTHOR = {Canuto, C. and Hussaini, M. Y. and Quarteroni, A. and Zang, T. A.},
     TITLE = {Spectral methods},
    SERIES = {Scientific Computation},
      NOTE = {Fundamentals in single domains},
 PUBLISHER = {Springer-Verlag, Berlin},
      YEAR = {2006},
     PAGES = {xxii+563},
      ISBN = {978-3-540-30725-9; 3-540-30725-7},
   MRCLASS = {65-02 (65M70 65N35 76M22)},
  MRNUMBER = {2223552},
MRREVIEWER = {R\'{e}mi Vaillancourt},
}

@book {Quarteroni2017,
    AUTHOR = {Quarteroni, Alfio},
     TITLE = {Numerical models for differential problems},
    SERIES = {MS\&A. Modeling, Simulation and Applications},
    VOLUME = {16},
      NOTE = {Third edition},
 PUBLISHER = {Springer, Cham},
      YEAR = {2017},
     PAGES = {xvii+681},
      ISBN = {978-3-319-49315-2; 978-3-319-49316-9},
   MRCLASS = {65-02 (65Mxx 65Nxx)},
  MRNUMBER = {3702005},
       DOI = {10.1007/978-3-319-49316-9},
}

@misc{kovachki2021neural,
      title={Neural Operator: Learning Maps Between Function Spaces}, 
      author={Nikola Kovachki and Zongyi Li and Burigede Liu and 
              Kamyar Azizzadenesheli and Kaushik Bhattacharya and Andrew Stuart and Anima Anandkumar},
      year={2021},
      eprint={2108.08481},
      archivePrefix={arXiv},
      primaryClass={cs.LG}
}

@book{Morrey,
    AUTHOR = {Morrey, Jr., Charles B.},
     TITLE = {Multiple integrals in the calculus of variations},
    SERIES = {Classics in Mathematics},
 PUBLISHER = {Springer-Verlag, Berlin},
      YEAR = {2008},
     PAGES = {x+506},
      ISBN = {978-3-540-69915-6},
   MRCLASS = {49-02 (49J10 49J45)},
  MRNUMBER = {2492985},
       DOI = {10.1007/978-3-540-69952-1},
}

@misc{lu2021comprehensive,
      title={A comprehensive and fair comparison of two neural operators (with practical extensions) based on FAIR data}, 
      author={Lu Lu and Xuhui Meng and Shengze Cai and Zhiping Mao and Somdatta Goswami and Zhongqiang Zhang and George Em Karniadakis},
      year={2021},
      eprint={2111.05512},
      archivePrefix={arXiv},
      primaryClass={physics.comp-ph}
}

@ARTICLE{Lanthaler2021,
    AUTHOR = {Lanthaler, Samuel and Mishra, Siddhartha and Karniadakis,
              George E.},
     TITLE = {Error estimates for {D}eep{ON}ets: a deep learning framework
              in infinite dimensions},
   JOURNAL = {Trans. Math. Appl.},
  FJOURNAL = {Transactions of Mathematics and its Applications. A Journal of
              the IMA},
    VOLUME = {6},
      YEAR = {2022},
    NUMBER = {1},
     PAGES = {tnac001, 141},
   MRCLASS = {47N10 (35Q68 47N30 68T07)},
  MRNUMBER = {4395087},
       DOI = {10.1093/imatrm/tnac001},
}

@misc{cai2021physicsinformed,
    AUTHOR = {Cai, Shengze and Mao, Zhiping and Wang, Zhicheng and Yin,
              Minglang and Karniadakis, George Em},
     TITLE = {Physics-informed neural networks ({PINN}s) for fluid
              mechanics: a review},
   JOURNAL = {Acta Mech. Sin.},
  FJOURNAL = {Acta Mechanica Sinica},
    VOLUME = {37},
      YEAR = {2021},
    NUMBER = {12},
     PAGES = {1727--1738},
      ISSN = {0567-7718},
   MRCLASS = {76-10},
  MRNUMBER = {4412280},
       DOI = {10.1007/s10409-021-01148-1},
       URL = {https://doi.org/10.1007/s10409-021-01148-1},
}

@article{Ruf2020,
    title={Neural networks for option pricing and hedging: a literature review},
    author={Johannes Ruf and Weiguan Wang},
        journal={J. Comp. Finance},
        year={2020},
        number={1},
        volume={24},
        doi={10.21314/JCF.2020.390},
        pages={1--46}
}

@article{DeepRitz,
    AUTHOR = {E, Weinan and Yu, Bing},
     TITLE = {The deep {R}itz method: a deep learning-based numerical
              algorithm for solving variational problems},
   JOURNAL = {Commun. Math. Stat.},
  FJOURNAL = {Communications in Mathematics and Statistics},
    VOLUME = {6},
      YEAR = {2018},
    NUMBER = {1},
     PAGES = {1--12},
      ISSN = {2194-6701},
   MRCLASS = {35J20 (35J91 65N30)},
  MRNUMBER = {3767958},
       DOI = {10.1007/s40304-018-0127-z},
}

@article{PiNNs,
    AUTHOR = {Raissi, M. and Perdikaris, P. and Karniadakis, G. E.},
     TITLE = {Physics-informed neural networks: a deep learning framework
              for solving forward and inverse problems involving nonlinear
              partial differential equations},
   JOURNAL = {J. Comput. Phys.},
  FJOURNAL = {Journal of Computational Physics},
    VOLUME = {378},
      YEAR = {2019},
     PAGES = {686--707},
      ISSN = {0021-9991},
   MRCLASS = {65M70 (68T05)},
  MRNUMBER = {3881695},
       DOI = {10.1016/j.jcp.2018.10.045},
}

@article{BeckEJentzen,
    AUTHOR = {Beck, Christian and E, Weinan and Jentzen, Arnulf},
     TITLE = {Machine learning approximation algorithms for high-dimensional
              fully nonlinear partial differential equations and
              second-order backward stochastic differential equations},
   JOURNAL = {J. Nonlinear Sci.},
  FJOURNAL = {Journal of Nonlinear Science},
    VOLUME = {29},
      YEAR = {2019},
    NUMBER = {4},
     PAGES = {1563--1619},
      ISSN = {0938-8974},
   MRCLASS = {65Z05 (35K55 35R60 60H10)},
  MRNUMBER = {3993178},
       DOI = {10.1007/s00332-018-9525-3},
}

@article{ChenChen1993,
author = {Chen, T. and Chen, H.},
title = {Approximations of continuous functionals by neural networks with application to dynamic systems},
journal = {{IEEE Transactions on Neural Networks}}, 
volume = {4},
year = {1993}, 
pages = {910 -- 918},
}

@article{MOPS20_938,
	author = {Marcati, Carlo and Opschoor, Joost A. A. and Petersen, Philipp C. and Schwab, Christoph},
	doi = {10.1007/s10208-022-09565-9},
	isbn = {1615-3383},
	fjournal = {Foundations of Computational Mathematics},
	journal = {Found. Comput. Math.},
	title = {Exponential {ReLU} Neural Network Approximation Rates for Point and Edge Singularities},
	url = {https://doi.org/10.1007/s10208-022-09565-9},
	year = {2022}}

@Article{OPS20_2738,
    author = {Joost A. A. Opschoor and Philipp C. Petersen and Christoph Schwab},
    title = {{Deep ReLU Networks and High-Order Finite Element Methods}},
    journal = {Anal. Appl. (Singap.)},
    volume = {18},
    number = {5},
    pages = {715--770},
    doi = {https://doi.org/10.1142/S0219530519410136},
    year = {2020}
}

@article{SchwabZech,
    AUTHOR = {Ch. Schwab and J. Zech},
     TITLE = {Deep learning in high dimension: neural network expression
              rates for generalized polynomial chaos expansions in {UQ}},
   JOURNAL = {Anal. Appl. (Singap.)},
  FJOURNAL = {Analysis and Applications},
    VOLUME = {17},
      YEAR = {2019},
    NUMBER = {1},
     PAGES = {19--55},
      ISSN = {0219-5305},
   MRCLASS = {68Q32 (35R60 41A25 41A46)},
  MRNUMBER = {3894732},
       DOI = {10.1142/S0219530518500203},
}

@article{EGPGSurvey,
    AUTHOR = {Elbr\"{a}chter, Dennis and Perekrestenko, Dmytro and Grohs,
              Philipp and B\"{o}lcskei, Helmut},
     TITLE = {Deep neural network approximation theory},
   JOURNAL = {IEEE Trans. Inform. Theory},
  FJOURNAL = {Institute of Electrical and Electronics Engineers.
              Transactions on Information Theory},
    VOLUME = {67},
      YEAR = {2021},
    NUMBER = {5},
     PAGES = {2581--2623},
      ISSN = {0018-9448},
   MRCLASS = {68T07},
  MRNUMBER = {4282380},
       DOI = {10.1109/TIT.2021.3062161},
}

@article{PETERSEN2018296,
title = "Optimal approximation of piecewise smooth functions using deep {ReLU} neural networks",
journal = "Neural Netw.",
volume = "108",
pages = "296 - 330",
year = "2018",
issn = "0893-6080",
doi = "10.1016/j.neunet.2018.08.019",
url = "http://www.sciencedirect.com/science/article/pii/S0893608018302454",
author = "Philipp Petersen and Felix Voigtlaender"
}

@unpublished{Costabel2010,
  TITLE = {{Corner Singularities and Analytic Regularity for Linear Elliptic Systems. Part I: Smooth domains.}},
  AUTHOR = {Costabel, Martin and Dauge, Monique and Nicaise, Serge},
  URL = {https://hal.archives-ouvertes.fr/hal-00453934},
  NOTE = {211 pages},
  HAL_LOCAL_REFERENCE = {v2: Improvement of layout},
  YEAR = {2010},
  MONTH = Feb,
  HAL_ID = {hal-00453934},
  HAL_VERSION = {v2},
}

@article{Opschoor2019,
  author = {J. A. A. Opschoor and Ch. Schwab and J. Zech},
  title = {Exponential {ReLU} {DNN} expression of holomorphic maps in high dimension},
  fjournal = {Constructive Approximation},
  journal = {Constr. Approx.},
  doi = {10.1007/s00365-021-09542-5},
  year = {2021}
}

@ARTICLE{Kutyniok2019,
       author = {{Kutyniok}, Gitta and {Petersen}, Philipp and {Raslan}, Mones and {Schneider}, Reinhold},
        title = "{A Theoretical Analysis of Deep Neural Networks and Parametric PDEs}",
      fjournal = {Constructive Approximation},
      journal = {Constr. Approx.},
         year = 2021,
        month = jun,
        doi = {10.1007/s00365-021-09551-4}
}

@article{lu2020deeponet,
Author = {Lu, Lu and Jin, Pengzhan and Pang, Guofei and Zhang, Zhongqiang and
   Karniadakis, George Em},
Title = {{Learning nonlinear operators via DeepONet based on the universal
   approximation theorem of operators}},
Journal = {Nature Machine Intelligence},
Year = {2021},
Volume = {3},
Number = {3},
Pages = {218--229},
Month = {{MAR}},
DOI = {10.1038/s42256-021-00302-5},
EISSN = {2522-5839},
}
\end{document}